\documentclass {article}

\usepackage{graphicx,amssymb}
\usepackage{amsmath}
\usepackage{amssymb}
\usepackage{amsthm}
\usepackage{amsfonts}
\newtheorem{corollary}{Corollary}
\newtheorem{theorem}[corollary]{Theorem}
\newtheorem{lemma}[corollary]{Lemma}
\newtheorem{conjecture}{Conjecture}
\newtheorem{proposition}[corollary]{Proposition}

\numberwithin{equation}{section}

\title{The conformally invariant measure on self-avoiding loops}
\author{Wendelin Werner}
\date {}

\begin {document}

\maketitle

\begin {abstract}
We show that there exists a unique (up to multiplication by constants) and natural measure 
on simple loops in the plane and on each Riemann surface, such that the measure is conformally invariant and 
also invariant under restriction (i.e. the measure on a Riemann surface $S'$ that is contained in 
another Riemann surface $S$, is just the measure on $S$ restricted to those loops that stay in $S'$).
We study some of its properties and consequences concerning
outer boundaries of critical percolation clusters and Brownian loops.
\end {abstract}

\font \m=msbm10
\newcommand{\R}{{\hbox {\m R}}}
\newcommand{\C}{{\hbox {\m C}}}
\newcommand{\Z}{{\hbox {\m Z}}}

\newcommand{\U}{{\hbox {\m U}}}
\def\H{{\hbox {\m H}}}

\def \eps {\varepsilon}
\def \cal {\mathcal}

\section {Introduction}

The aim of the present paper is to construct and describe a natural measure 
on the set of self-avoiding loops in the plane and on any Riemann surface. By a self-avoiding loop on a 
surface $S$, we mean a continuous injective map from the unit circle into $S$ modulo monotone reparametrizations (i.e. we look only at the trace of the loop and forget about its parametrization).
We will construct a measure that possess some strong conformal invariance
 properties, and we shall see that this measure is the only one with these
 properties. 

Let us first describe this strong conformal invariance property, phrased in terms of a measure on the set of 
self-avoiding loops in the plane: We say that such a measure $\mu$ satisfies conformal restriction if for any two 
conformally equivalent domains $D$ and $D'$ (i.e. such that there exists a conformal map from 
$D$ onto $D'$) in the plane, the image of the measure $\mu$ restricted to the set of loops that stay in $D$, via any conformal map $\Phi$ from $D$ onto $D'$, is exactly the measure $\mu$ restricted to the set of loops that stay in $D'$.

Note that this condition implies in particular that the measure $\mu$ is translation-invariant, scale-invariant (and therefore that it has infinite total mass provided that $\mu \not= 0$), and that if $\mu$ satisfies conformal restriction, then so does $c\mu$ for any positive constant $c$. 
We are going to say that a measure on the set of self-avoiding loops is non-trivial if for some $0< \delta<\Delta < \infty$, the mass of the set of loops of diameter at least $\delta$ and that stay in some disc of radius $\Delta$ is neither $0$ nor infinite.
 
We shall prove the following result:

\begin {theorem}
\label {T0}
Up to multiplication by a positive constant, there exists a unique non-trivial measure on the set of self-avoiding loops in the plane that satisfies conformal restriction.
\end {theorem} 

Before giving some motivation, let us now describe the counterpart of this result in 
term of measures on loops on Riemann surfaces.
Suppose that for each Riemann surface $S$ (we do not require $S$ to be closed), we are given a measure $\mu_S$ on the set of self-avoiding loops on $S$. We say that the family of  measures $(\mu_S)$ satisfies conformal restriction if the following two conditions hold:
\begin {itemize}
\item
For any conformal map $\Phi$ from one surface $S$ onto another surface $S'$, 
$\Phi \circ \mu_S = \mu_{S'}$
\item
For any $S$ and any $S' \subset S$, $\mu_{S'}$ is equal to the measure $\mu_S$ restricted to those loops 
that stay in $S'$.
\end {itemize}

Suppose that the family $(\mu_S)$ 
satisfies conformal restriction. It is then clear that the measure $\mu=\mu_\C$  in this family satisfies conformal restriction in the previous sense. 
It is also easy to see that for any Riemann surface $S$, one can define a countable family of annular regions (i.e. conformally equivalent to a planar annular region) such that any self-avoiding loop in $S$ is contained in at least one of these regions. Hence, the knowledge of the $\mu_D$'s for all annular regions $D$ determines the entire family $(\mu_S)$. One consequence is that the entire family $(\mu_S)$ is described by the measure
$\mu=\mu_\C$. This will lead to:

\begin {theorem}
\label {T1}
Up to multiplication by a positive constant, there exists a unique family of such measures $(\mu_S)$ that satisfies conformal restriction. 
\end {theorem}

\medbreak

The proof of these results will consist of two different parts. First, one shows that there exists at most one measure $\mu$ that satisfies conformal restriction, and then one constructs explicitely a measure $\mu$ that possesses this property.
As we shall see, the first step is not very difficult and requires no 
prerequisite. On the other hand, the existence part of the proof of the theorem is non-trivial.
Recall that we are not restricting ourselves to simply connected domains $D$; for instance, the measure
$\mu$ will be invariant under the inversion $z \to 1/z$ from $\C \setminus \{ 0 \}$ onto itself.

\begin {figure}
\label {f1}
\begin {center}
\includegraphics [width=2.5in]{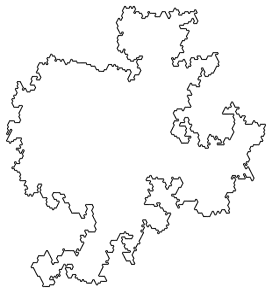}
\caption {A self-avoiding loop (sketch)}
\end {center}
\end {figure}

\medbreak

Let us now list some other properties of the measures described in Theorems \ref {T0} and \ref {T1} that we shall 
derive in this paper:
\begin {itemize}
\item
For any given simply connected domain $D \subset \C$, for any simply connected $D' \subset  D$, and $z \in D'$, 
consider the conformal map $\Phi$ from $D'$ onto $D$ that fixes $z$ and that has a  positive derivative at $z$.
Then 
$$ \mu_D ( \{ \gamma \ : \ \gamma \hbox { disconnect z from } \partial D
\hbox { and } \gamma \notin D' \} ) = c \log \Phi'(z)$$
for some $c$. This formula (for fixed $D$ and letting $z$ and $D'$ vary)  characterizes the measure $\mu_D$ fully (and therefore $\mu$ too because of its scale-invariance). 

\item The measure $\mu$ can be viewed as the measure on outer boundaries of planar Brownian loops in the plane, where these Brownian loops are defined under the ``Brownian loop measure'' introduced in \cite {LW}. This is roughly speaking a scale-invariant and translation-invariant measure on (unrooted) Brownian loops in the plane.

\item
For each Riemann surface $S$ (and in particular for $S=\C$), the measure $\mu_S$ is supported on the set of loops with Hausdorff dimension $4/3$ (see Figure \ref {f1}). This will follow from the construction of $\mu$ via outer boundaries of Brownian loops that we have just mentioned (recall that it has been proved-- see
 \cite{LSW2,LSW4/3}-- that the dimension of outer boundaries of two-dimensional
 Brownian motions is $4/3$).

\item One can also view the measures $\mu_S$ as measures on ``Schramm-Loewner Evolutions loops'' with parameter $8/3$. This is related to the special properties of SLE with that parameter (e.g. \cite {LSWr}). It will be an instrumental fact in the proof of Theorems \ref {T0} and \ref {T1}.

\item
One can also interpret the measure $\mu$ as the scaling limit of the measure on critical percolation cluster outer perimeters. In particular, this will show that the shapes of percolation cluster outer perimeters in the scaling limit are exactly the same (in law) as that of Brownian loops. 

\item    
The $\mu$ mass of the set of loops that stay and go ``around'' the annulus
$\{ z \ : \ |z| \in (1, \exp (\rho))\}$ decays like a constant times
$ \exp \{-  5 \pi^2 / 4\rho \}$ when $\rho \to 0+$.
This estimate also holds for other annular regions by conformal invariance. 
\end {itemize}

This last estimate can be related to the conjecture that the measures $\mu_S$ are the scaling limits of simple measures on self-avoiding polygons on fine-mesh discrete approximations of the surface $S$ (see section \ref {section6.1}  and \cite {LSWsaw}). Roughly speaking, $c\mu$ for some $c$ should be the fine-mesh limit (i.e. $\delta \to 0$) of the measure that assign a mass $\lambda^{-n(l)}$ to each self-avoiding loop $l$ on the grid $\delta \Z^2$, where $n(l)$ denotes the number of steps of the loop $l$, and $\lambda$ is the connectivity constant of the lattice $\Z^2$.

The description of measures $\mu$ in terms of Brownian outer boundaries or percolation 
outer perimeters in the scaling limit, combined with the properties of the measures $\mu$ (such for instance as the 
invariance under inversion $z \mapsto 1/z$) has in turn rather surprising consequences for Brownian loops and shapes of percolation clusters. In particular, it will show that the shapes of
 ``outer'' and ``inner'' boundaries have exactly the same law.

The paper will be structured as follows: in the next section, which can be viewed as the second part of the introduction,  we make some further general comments and recall a few general facts and ideas from the papers \cite {LW} and \cite{LSWr}.  
In Section \ref {s3}, we show uniqueness of the measure $\mu$ satisfying a weaker form of conformal restriction  (namely  that the restrictions of $\mu$ to simply connected sets are all conformally equivalent).
This implies in particular the uniqueness part of Theorem \ref {T0}.
 In the subsequent section, 
 we will then show how to construct a measure $\mu$ satisfying this weak conformal restriction property
 using the Brownian loop
measures (to each Brownian loop we associate its outer boundary which is a self-avoiding loop). We will then focus on the restriction of this measure to annular regions. In particular, 
using SLE$_{8/3}$ considerations, we shall prove that for any two conformally equivalent annular regions in the plane, the restriction of $\mu$ to the second one is equal to the conformal image of the restriction of $\mu$ to the first one.
This  will in turn imply Theorem \ref {T0} and show that it is 
possible to define $\mu_S$ for any surface in such a way that the obtained family of measures satisfies conformal restriction (i.e. that Theorem \ref {T1} holds too).
In the final two sections, we will study asymptotics of the mass of loops that go around thin annuli (motivated by the 
discrete self-avoiding walks problems) and the relation to critical percolation.

\section {Further motivation and background}

This paper builds on ideas of the joint paper with Greg Lawler and Oded Schramm on conformal restriction in the chordal case \cite {LSWr}, and uses also the Brownian loop measure introduced in joint work with Greg Lawler \cite {LW}. 
 A survey of the results concerning SLE with parameter $8/3$ and Brownian measures (on excursions and loops) that we will use can for instance be found in the first four sections of \cite {Wcrr}. 
 Let us however briefly recall some aspects of the two papers
 \cite {LW,LSWr} now in order to put our results into perspective:
 
 \medbreak
 {\bf Self-avoiding loops and Brownian loops.}
Because of its conformal invariance,
planar Brownian motion can be a useful tool to interpret and prove results 
in complex analysis (e.g. \cite {Bass}) even though many questions in this field 
were already settled long before Brownian motion was
properly defined. Usually, by conformal invariance of planar Brownian motion one means  that 
the conformal image of a Brownian motion started at a given point in a given domain is a (time-changed) Brownian motion in the image domain started
at the image of the starting point (this had been first observed by Paul L\'evy in the 1950's and it can be viewed as a
direct consequence of It\^o's formula, see e.g. \cite {LGln}).   
If one does however not prescribe any starting point on a Riemann surface and still wants to state a conformal invariance property for measures on Brownian paths,
it is natural to look for measures on Brownian loops. 
In \cite {LW}, we defined a measure on Brownian loops in the complex plane (and 
subsets of the complex plane) that is indeed conformally invariant. More precisely, it satisfies a property similar to the conformal restriction that we described before (but of course, Brownian loops are not simple loops): 
The image under a conformal map $\Phi$ of the measure on Brownian loops restricted to those that stay in a domain $D \subset \C$ is exactly the measure on Brownian loops 
restricted to those loops that stay in $\Phi (D)$. 
As we shall see later in this paper, 
it is  easy to generalize the definition of the Brownian loop measure to 
 general Riemann surfaces and to show that this conformal restriction 
property still holds. 

In the case of hyperbolic surfaces (where Green's functions are finite) such for instance as the bounded subsets of $\C$, 
quantities like the mass of Brownian loops that intersect two disjoint sets turn out to be 
rather natural conformal invariants that are
for instance related to Schwarzian derivatives (see e.g. \cite {LSWr})
and they show up under various guises in 
conformal field theory. Furthermore \cite {LW,Wls}, it is natural to construct Poissonized samples of 
these Brownian loop-measures (that we call the Brownian loop-soups in \cite {LW}) 
and to study the geometry of the obtained sets \cite {Wls}. This gives a way to construct the so-called conformal loop ensembles (CLE) that are conjectured to be the scaling limits of various critical two-dimensional models from statistical physics, see  \cite {ShW,Wln3}.

The Brownian loop-measure on compact surfaces (for which Green's functions are infinite),  such for instance the Riemann sphere,  
is still easy to define but not so easy to work with. 
For instance, the mass of the set of loops that intersect two disjoint open sets is 
infinite because of the too many long loops (due to the recurrence of the Brownian motion).  
This implies in particular that the corresponding Brownian loop-soups become rather uninteresting, because there always exists just one dense cluster of 
loops. However, it should be possible to define the CLEs and conformal field theories on such surfaces in order to describe scaling limits of lattice models.

We shall see that the measure on self-avoiding loops $\mu_S$ that we are defining is in a way
better-suited to higher-genus surfaces and compact surfaces than the Brownian loop-measure. 
For simply connected planar domains, the self-avoiding loop measure $\mu$ is very directly related to the Brownian loop
measure (it is the measure on ``outer boundaries of the Brownian loops'') and 
therefore also to Schwarzian derivatives etc. Its Poissonian samples also define the same CLEs as the Brownian loop-measure does.
But these nice properties will still hold for any Riemann surface. For instance, it 
allows to construct directly loop-clusters for a loop-soup in the entire plane.
In this setup, all the nice properties of the obtained
cluster boundaries are just consequences of the conformal restriction property of 
the family $(\mu_S)$.

\medbreak
{\bf Chordal restriction.}
Let us now briefly recall some aspects of conformal restriction in the chordal case, as 
studied in \cite {LSWr}. There, we focused on probability measures on self-avoiding curves that join two 
prescribed boundary points of simply connected domains (we will only discuss simply connected domains in this paragraph). The law of such a curve was therefore characterized by the 
domain $D$ and also by the two end-points $A$ and $B$ of the path. We showed that there exists a unique
family $(P_{D,A,B})$ of  
probability measures on that set of self-avoiding curves, that satisfies conformal invariance
(i.e. $\Phi \circ P_{D,A,B} = P_{\Phi(D), \Phi (A), \Phi(B)}$) and restriction in the sense that the law
of the curve in $D$, {\em conditioned} to stay in $D'$ was equal to the law in the smaller domain with the boundary points $A$ and $B$.  The constraint on $D'$ is that $A$ and $B$ are also boundary points of $D'$ and for instance that 
$D \setminus D'$ is at positive distance of $A$ and $B$ (to ensure that the curve stays in $D'$ with positive probability so that the conditioning makes sense). 
Furthermore (see \cite {LSWr} for more details), we showed that
$$ P_{D, A, B} ( \gamma \subset D' ) = (\Phi' (A) \Phi'(B))^{5/8} 
$$
when $\Phi$ is a conformal map from  $D' \subset D$ back onto $D$ that preserves the two boundary points $A$ and $B$.

One can then 
try to define a family of measures $\mu_{D,A,B}$ on such curves in such a way that $\mu_{D,A,B}$ restricted to 
those curves that stay in $D'$ is exactly equal to $\mu_{D',A,B}$. In other words, we allow the measure to have a total mass different to 1, but we require that the restriction property holds with no renormalizing term (for the family $P_{D,A,B}$, there is such a term due to the conditioning).
Let us add the technical condition  that we restrict the definition of $\mu_{D,A,B}$ to 
those sets $D$ and points $A$, $B$, such that the boundary of $D$ is smooth in the neighborhood of $A$ and $B$.
Then, we simply define 
$$\mu_{D,A,B}= \Phi' (A)^{-5/8} \Phi' (B)^{-5/8} P_{D,A,B},$$
 where $\Phi$ is a conformal map
from some fixed reference domain (say $\H, 0, \infty$) onto $D,A,B$. The family of measures $\mu_{D,A,B}$ is then defined modulo a multiplicative constant (corresponding to the choice of the reference domain) and can be shown to be the unique one (on the set of simple paths from one boundary point to another of smooth simply connected domains) satisfying this restriction property.

It was very natural to study this chordal restriction property in the light of the definition and properties of 
the chordal Schramm-Loewner evolutions that are precisely random paths joining one boundary point of a domain to 
another. In fact, the probability $P_{D,A,B}$ is 
exactly  the law of the SLE with parameter $8/3$ in $D$ from $A$ to $B$.  It also turned  out (see e.g. \cite {FW,BB}) that this chordal restriction property of SLE with parameter $8/3$ had interpretations in terms of boundary conformal field theory, and that it yields a  simple definition of
 the SLE$_{8/3}$ in non-simply connected domains.
Also, a by-product of the results of \cite {LSWr} was the fact that the Brownian outer boundary looked ``locally'' like that of a critical percolation boundary, and also like a SLE$_{8/3}$, because of global identities in law between outer boundaries of the union of 5 Brownian excursions and that of the union of 8 SLE$_{8/3}$'s. This showed that ``locally'', the inside boundaries and the outside boundaries of the Brownian loops and of percolation clusters looked like an SLE$_{8/3}$. 

So, we see that part of the present paper is really a continuation of some of the ideas developed in \cite {LSWr}.
In a way,
 just as for the Brownian loops compared to Brownian motions, the measure
 on self-avoiding 
loops $\mu$ in the plane is invariant under a wider class of transformations than the chordal measures, because no marked points are prescribed. This will make it well-suited to general Riemann surfaces and to the description and construction of conformally invariant quantities and objects. 

\medbreak

{\bf Further motivation.}
It is also worth emphasizing that the conformal restriction property enables us to 
understand the behavior of these measures when one distorts locally the Riemann surface (for instance, one cuts out a small disc and glues some other structure), so that $\mu$ provides 
conformally invariant information that behaves nicely under perturbation of the conformal structure.  
We  believe that it is a natural and hopefully fruitful way to describe useful conformal invariant quantities in this setup \cite {Wip}.

Also, one can note (see e.g. \cite {KY}) that the space of shapes of self-avoiding loops
 (as we shall define it later in the paper) 
can be realized in many different ways in particular involving the Lie algebra Vect$(S^1)$, 
 and that it has been the subject of many investigations 
(see e.g. \cite {Mall}).
We hope to clarify the link  with the present paper in \cite {Wip}.
Measures on loops
 should also be the natural setup for considerations related to the geometric-theoretical
aspects of SLE as looked into in  \cite {FK,F,Ko}.

\section {Simply connected case: Uniqueness}
\label {s3}

This section requires no prerequisite (i.e. no knowledge on Brownian loops nor on SLE).

Let us settle a few technical details:
Throughout the paper, when we will talk about self-avoiding loops, we will not care about their time-parametrizations.
In other words, we identify two self-avoiding loops as soon as their traces are the same. This means that two injective continuous maps $\gamma^1$ and $\gamma^2$ from the unit circle $S^1$ into some domain define the same loop if there exists a continuous monotone bijection $w$ from $S^1$ onto $S^1$ such that $\gamma^1=\gamma^2 \circ w$. 
Note that with this definition, a path of zero length (i.e. ``a point'') is not a self-avoiding loop. 

This set of self-avoiding loops can be endowed with various natural metrics. We will not discuss this here, but we have at least to specify precisely what $\sigma$-field we will implicitly refer to when we define measures on this set. 
 Throughout the paper, we will say that a bounded open set $A$ is an annular region (not to be confused with exact annuli) in the plane if it is conformally equivalent to some annulus $\{z \ : 1 < |z| < R \}$. In other words, $A$ is a bounded open set and $\C \setminus A$ has two connected components (that are not singletons). The bounded connected component of $\C \setminus A$ will be called the ``hole'' of $A$.
For each such annular region $A$, define the set ${\cal U}_A$ of self-avoiding loops that stay in $A$ and that have a non-zero index around the hole of $A$ (or in other words, that disconnects the two connected components of $\partial A$). We will denote by 
${\cal G}$ the $\sigma$-field generated by this family. This is the $\sigma$-field that we will work with (we leave it to the interested reader to check the relation between ${\cal G}$ and the Borel $\sigma$-field induced by the Hausdorff distance on the set of self-avoiding loops).

We will also use the countable family ${\cal A}$ of annular regions $A$ such that the inner and the outer boundaries of $A$ consist of finite loops that are drawn on the lattice $2^{-p} \Z^2$ for some $p \ge 1$. Clearly, since any annular region can be seen as the increasing limit of such dyadic annular regions, the family $({\cal U}_A, A \in {\cal A})$ generates ${\cal G}$ too. The sets of self-avoiding loops that we shall consider can be easily shown to be in ${\cal G}$ (for instance because they are countable unions of  ${\cal U}_A$'s).

When $\mu$ is a measure supported 
on the set of self-avoiding loops in the complex plane, then, for each simply connected domain $D$ (unless otherwise stated, domains are open sets), we simply
define $\mu_D$ to be the measure $\mu$ restricted to the set of 
loops that stay in $D$. Throughout this paper, by a slight abuse of notation, we will write
$\mu_D = \mu 1_{\{ \gamma \subset D\}}$ (instead of $\mu_D ( \cdot) = \mu ( \cdot \cap \{ \gamma \subset D \})$).

Recall that we say that a measure on the set of self-avoiding loops is non-trivial if for some $0< \delta < \Delta < \infty$, the mass of the set of loops of diameter at least $\delta$ that are contained some disc of radius $\Delta$, is neither zero nor infinity. 

We will also use a weaker form of conformal restriction in this paper:
We will say that a measure $\mu$  on the set of self-avoiding loops in the plane satisfies {\em weak}
conformal restriction if for any two
simply connected 
domains $D$, $D'$ and any conformal map $\Phi$ from $D$ onto $D'$,
$$ \Phi \circ \mu_D  = \mu_{D'}.$$
The difference with the (strong) conformal restriction is that here, the domains $D$ and $D'$ are simply connected.

We can now state the main result of this section:

\begin {proposition}
\label {P1}
Up to multiplication by a positive constant, there exists at most one non-trivial measure $\mu$ (on the set of self-avoiding loops) satisfying weak conformal restriction.

Furthermore, for any simply connected sets $\tilde D \subset D$, and any $z \in \tilde D$,
\begin {equation}
 \mu ( \{ \gamma \ : \ \gamma \subset D, \ \gamma \not\subset \tilde D , \ 
\gamma \hbox { disconnects z from }\partial D  \} )
= c \log \Phi' (z),
\label {phi'}
\end {equation}
where $c$ is a positive constant (depending on $\mu$ only), and $\Phi$ denotes the conformal map 
from $\tilde D$ onto $D$ such that $\Phi(z)=z$ and $\Phi'(z)$ is a positive real.
\end {proposition}

This result does not say anything about whether this measure exists.

\medbreak

We will also prove the following closely related  lemma. 
It is a similar result for measures on loops that surround a given point.
Suppose that $\nu$ is a non-trivial measure supported on the set of 
self-avoiding loops in the complex plane that disconnect $0$ from infinity.
For any simply connected domain $D$ containing the origin, we define 
again $\nu_D = \nu 1_{\gamma \subset D}$. 

\begin {lemma}
\label {L1}
Up to multiplication by a positive constant, there exists at most one such measure $\nu$ such 
 that for any simply connected domain $D$ containing the origin and any
conformal map $\Phi$ defined on $D$ such that $\Phi(0) = 0$, one has
$\Phi \circ \nu_D = \nu_{\Phi (D)}$.
Furthermore, there exists a positive constant $c$ such that 
for any simply connected subset $U$ of the unit disc $\U$ such that 
$0 \in U$, 
\begin {equation}
\label {E4} \nu (\{\gamma \ : \ \gamma \subset \U , \gamma \not\subset U\} ) 
= c \log (\Phi' (0))
\end {equation}
where $\Phi$ is the conformal map from $U$ onto $\U$ such that $\Phi(0)=0$
and $\Phi'(0)>0$.
\end {lemma}

\medbreak
\noindent
We will first prove the lemma.

\medbreak
\noindent
{\bf Proof of the Lemma.}
Suppose that $\nu$ is a measure supported on the set of self-avoiding loops that surround the origin,
such that for any simply connected domain $D$ containing the origin and any conformal map 
$\Phi : D \to \Phi (D)$ with $\Phi( 0)= 0$, one has $\nu_{\Phi (D)} = \Phi \circ \nu_D $.
Our first goal is to prove that it necessarily satisfies (\ref {E4}).

Consider a simply connected subset $U$ of the unit disc that contains the origin.
We define the unique conformal map $\varphi=\varphi_U$ from $U$ onto $\U$ with 
$\varphi (0)= 0 $ and $\varphi' (0) > 0$. Note that as $U \subset \U$, 
$\varphi'(0)  \ge 1$ (recall for instance that if $Z$ is a planar Brownian motion 
started from the origin and stopped at its first exit time $T$ of $U$,
then $\log (\varphi' (0)) = - E (\log | Z_T | ) \ge 0$). 

Note that $\varphi$ describes $U$ fully as $U = \varphi^{-1} ( \U)$. 
We now define
$$ A (\varphi_U) = \nu (\{ \gamma \ : \ \gamma \subset \U, \ \gamma \not\subset U\} ).
$$

\begin{figure}[htbp]
\begin{center}
\includegraphics [height=1in]{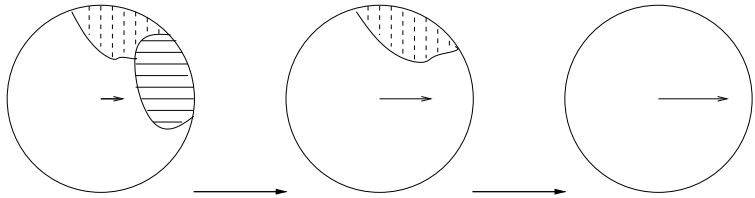}
\caption{$\varphi_V \circ \varphi_U$}
\label{figure:example1}
\end{center}
\end{figure}

This quantity measures the mass of the set of
loops  that stay in the unit disc but do exit the smaller domain $U$
(recall also that $\nu$ is supported on the set of loops that surround the origin).
Now, under the conditions of the lemma,
it is immediate to see that for any simply connected subsets $U$ and $V$ of the unit disc
 that do contain the 
origin,
\begin {equation}
\label {+}
 A ( \varphi_V \circ \varphi_U ) = A ( \varphi_U) + A (\varphi_V) .
\end {equation}
Indeed, 
$A (\varphi_V \circ \varphi_U)$ is the mass of loops that surround the origin, stay in $\U$, but hit 
$\U \setminus U$ or $\varphi_U^{-1} ( \U \setminus V)$. Note that these two sets are disjoint.
But by definition
\begin {equation}
\label {+1}
\nu (\{ \gamma \ : \ \gamma \subset \U, \ \gamma \not\subset U
 \}) = A (\varphi_U)
\end {equation}
and 
$$ 
\nu (\{ \gamma \ : \ \gamma \subset \U, \ \gamma \not\subset V
 \} )
= A (\varphi_V).
$$
By conformal invariance (via $\varphi_U^{-1})$, we see that this second quantity is also equal 
to
\begin {equation}
\label {+2}
\nu (\{ \gamma \ : \ \gamma \subset U, \ \gamma \not\subset \varphi_U^{-1} ( V)\} ) 
\end {equation}
Summing (\ref {+1}) and (\ref {+2}), we get (\ref {+}).

Let us now focus on a simple class of domains $U_t$ in $\U$: For each positive $t$, we let 
$U_t = \U \setminus [r_t, 1)$, where $r_t$ is the positive real such that 
$\varphi_{U_t}' (0) = e^t$ (it is not difficult to work out the explicit expression of 
$r_t$ but we do not need it here). 
Note that with this definition, the family $(\varphi_{U_t})_{t \ge 0}$ is a semi-group 
(i.e. $\varphi_{U_t} \circ \varphi_{U_s} = \varphi_{U_{s+t}}$)
because (by symmetry), the composition of these 
two maps is still a map from some $U_r$ onto $\U$, and $r = s+t$ because 
$(\varphi_{U_t} \circ \varphi_{U_s})' (0) = e^{s+t}$.
Hence, $t \mapsto  A (\varphi_{U_t})$ is a non-decreasing function 
from $(0, \infty)$ into $(0, \infty)$, such that 
$$ A (\varphi_{U_{t+s}}) = A (\varphi_{U_t}) + A(\varphi_{U_s}) .$$
This implies that for some positive constant $c$, we have that $A(\varphi_{U_t}) =c t
=c \log \varphi_{U_t}' (0)$ (the constant is positive and finite because $\nu$ is non-trivial).
In the sequel, we fix $c$.

For each $\theta \in [0, 2 \pi)$ and $t >0$, define the set
$$ U_{t, \theta}  = \U \setminus [ r_t e^{i \theta}, e^{i \theta}) .$$
Clearly, $\varphi_{U_{t, \theta}}' (0) = t$  and because of 
the invariance of $\nu$ under rotations, we have $A (\varphi_{U_{t, \theta}}) = A (
\varphi_{U_t}) = c t $.

Let ${\cal S}$ denote the semi-group of conformal maps generated by the family 
$(\varphi_{U_{t, \theta}}$, $t >0$, $\theta \in [0, 2 \pi))$ (i.e. the set of finite compositions 
of such maps).
Because of (\ref {+}), we deduce from the above that for any $\varphi \in {\cal S}$,
\begin {equation}
A ( \varphi) = c \log ( \varphi'(0)).
\label {phi'2}
\end {equation}

The family ${\cal S}$ is ``dense'' in the class of conformal maps $\varphi_U$ from some simply connected subset 
$U$ onto $\U$ in the following sense (this is a standard fact from the theory developed by Loewner):
For any $U$, there exists an increasing family $U_n$ in ${\cal S}$ such that 
$\cup_n U_n = U$. Clearly, it follows that $\varphi_{U_n}' (0) \to \varphi_U'(0)$.
Furthermore, because a loop is a compact subset of the complex plane, $\gamma$ exits $U$ if and only if it 
exits all $U_n$'s, so that 
$ A (\varphi_{U_n}) \to A ( \varphi_U)$ as $n \to \infty$.
Combining the above, we get that for any $U$,
$$
A ( \varphi_U ) = \lim_n A (\varphi_{U_n}) = c  \lim_n \log \varphi_{U_n}' (0) 
= c \log \varphi_U' (0).
$$
This shows that $\nu$ must indeed satisfy (\ref {E4}). 

\medbreak
 The lemma will then follow from 
 rather simple measure-theoretical considerations that we now detail. 
 Recall that a finite measure on a measurable space is characterized by its total mass and its value on a subset of the $\sigma$-field that generates the $\sigma$-field and that is stable under finite intersections. 
 
 Until the end of this proof, we are going to restrict ourselves to the set of loops that surround the origin. 
 
 We say that an annular region is ``good'' if the origin is in the interior of its inner hole. Note that the set of loops in ${\cal U}_A$ that do surround the origin (when $A$ is not good), can be viewed as a countable 
 unions of sets  ${\cal U}_{A_j}$, where the $A_j$s are good annular regions. 
 Hence, we can restrict ourselves to the $\sigma$-field generated by good annular regions. 
 
But the set of all ${\cal U}_A$, where $A$'s are good annular regions is stable under finite intersection. So, we would like to prove that the values $\nu ({\cal U}_A)$ for good annular regions are determined uniquely (once the constant $c$ has been chosen). Let us fix a good annular region $A$.

Conformal invariance implies that (under the previous assumptions for $\nu$) there exists $c>0$ such
that for any two simply connected domains $V$ and $U$ with $0 \in V \subset U \not= \C$,
the $\nu$ mass of the set of loops 
$$X (V, U) =  \{ \gamma \ : \gamma \hbox { surrounds } 0, \ \gamma \not\subset V, \ \gamma \subset U \}$$
is equal to $c \log \Phi' (0)$ where $\Phi$ is the conformal map from $V$ onto $U$ that fixes the origin such that $\Phi'(0)$ is a positive real (this follows immediately from (\ref {E4}) using the conformal equivalence of $\nu_U$ and $\nu_\U$).

We now construct another family of observables that is stable under finite intersections:
We now define $U_0$ to be the union of $A$ with its hole. We choose $V_0$ to be a simply connected set containing the origin that is at positive distance of $A$. 
Define $X_0 = X (V_0, U_0)$ and the family $\Pi_0=\Pi (V_0,U_0)$  of events $X(V_0,W) \subset X_0$ 
where $W$ spans all possible simply connected sets with $V_0 \subset W \subset U_0$.
Note that $\Pi_0$  is indeed stable under (finite) intersections. (\ref {E4}) for a given $c$ 
yields the value of $\nu$ on the family $\Pi_0$, and its total mass on $X_0$, and it therefore 
characterizes $\nu$ on the 
entire $\sigma$-field $\Sigma_0 = \Sigma(V_0,U_0)$ generated by $\Pi_0$. 

Let us denote the outer boundary of $A$ by $a_2=a_2(A)$ and its inner boundary by $a_1=a_1(A)$. 
Note that ${\cal U}_A \subset X_0$. 
It is easy to check that a self-avoiding loop is in ${\cal U}_A$ if and only if, for some $m \ge 1$ and all $p\ge 1$, it intersects any (of the countably many) simple continuous path on $2^{-p} \Z^2$ that intersects $a_2$ and that is at distance less than $2^{-m}$ from $a_1$. Hence ${\cal U}_A$ is a countable union of countable intersections  of sets $X(V_0, U_0 \setminus \eta)$, so that ${\cal U}_A \in \Sigma_0$. Hence, 
$\nu ({\cal U}_A)$ is in fact uniquely determined.

The family ${\cal U}_A$ when $A$ spans all good annular regions is also stable under finite intersections. The value of 
$\nu$ on this family is determined. Furthermore the masses $\nu({\cal U}_{A_0})$'s are finite for good annular regions $A_0$. It follows that for each good annulus $A_0$, $\nu$ restricted to ${\cal U}_{A_0}$ is fully determined. Since this is true for any $A_0$, it follows that $\nu$ itself is uniquely characterized once the constant $c$ has been chosen.
\qed

\medbreak
{\bf Proof of the proposition.}
Let us now assume that $\mu$ is a non-trivial measure on the set of self-avoiding loops in the plane that satisfies weak conformal restriction. Define $\nu$ to be its restriction to the set of loops that surround the origin. It does then  satisfies the conditions of the lemma. It is therefore uniquely determined and satisfies (\ref {E4}). By translation-invariance, this characterizes $\mu$ too (because for any self-avoiding loop $\gamma$, there exists a point $z$ with rational coordinates such that $\gamma -z$ surrounds the origin), and the formula (\ref {phi'}) follows from the expression of $\nu (X(V,U))$ and translation invariance.
\qed

\section {Simply connected case: Construction}

We have reduced the class of measures on loops that satisfy weak conformal restriction to a possibly
empty one-dimensional (i.e. via multiplication by  constants) family of measures. 
It is of course natural to try to construct such a measure, to make sure that this family is not empty.
We are going to do this using the Brownian loop-measure introduced with Greg Lawler in \cite {LW}.

\subsection {The Brownian loop-measure}

Since it is an instrumental tool throughout the present paper, it 
is worthwhile to recall the construction of the Brownian loop-measure defined in \cite {LW} (in $\C$ and 
subsets of $\C$) and its basic properties. 

It is well-known since Paul L\'evy that planar Brownian motion is conformally invariant. This is usually stated as 
follows. Suppose that $Z$ is a Brownian motion started from $z \in D$, and stopped at its first exit time
$T$ of 
the domain open $D \subset \C$.
 Suppose that $\Phi$ is a conformal map from $D$ onto $\Phi (D)$. Then, up to time-reparametrization, 
$\Phi ( Z [0,T])$ is a Brownian motion started from $\Phi(z)$ and stopped at its first exit time  of $\Phi(D)$.

Also, recall that $\log |Z_t|$ is a local martingale when $Z$ is a planar Brownian motion
started away from the origin.
It follows easily from these two facts that for each $z$, one can define a measure $N^z$ supported 
on Brownian loops that start and end at $z$ as follows:
Consider, for each $\eps > 0$, the law $P_{z, \eps}$ of a Brownian motion started
uniformly on the circle of radius $\eps$ around $z$ and stopped at its first hitting 
of the circle of radius $\eps^2$.
Then, consider the limit $N^z$ when $\eps \to 0$ of 
$$4 \log (1/\eps) \times P_{z, \eps}.$$
This limit is an infinite measure supported on the set of Brownian paths $Z$ that start and end at $z$.
The mass of the set of loops that reach the unit circle, but not the circle of radius $\exp r$ (for $r >0$) 
is $r$.
Conformal invariance of planar Brownian motion implies readily that this measure $N^z$ is also conformally invariant 
in the following sense. If we consider a domain $ D$ with $z \in D$ and a conformal map $\Phi : D \to \Phi (D)$, 
then 
\begin {equation}
\label {bm1} 
\Phi ( N^z 1_{ Z \subset D} ) = N^{\Phi(z)} 1_{Z \subset \Phi(D)}
\end {equation}
Here, as in the rest of this paper, we do not care about the time-parametrization of the Brownian loop. We identify two Brownian loops if one can be obtained from the other by an ``increasing'' reparametrization (which is 
more than just identifying the traces since the Brownian motions have double points -- as opposed to the self-avoiding loops).
However, and we shall use this in a moment, a Brownian loop $Z$ carries its ``natural'' time-parametrization given by its quadratic variation. It has a natural time-length $T(Z)$ which 
is not a conformally invariant quantity. In fact, it is easy to see that 
$$ T ( \Phi(Z)) = \int_0^{T(Z)} | \Phi' (Z_t)|^2 dt $$
when $Z$ is parametrized by its ``natural'' time.

Now, consider the product measure $\hat M = d^2 z \otimes N^0$, where as in the rest of this paper, $d^2 z$ denotes the 
Lebesgue measure in the complex plane.
This defines a couple $(z, Z^0)$, where $Z^0$ is a Brownian loop rooted at the origin, and $z \in \C$.
We then simply define $Z^z = z + Z^0$ which is now a Brownian loop rooted at $z$. Finally, we call $Z$ the equivalence class of $Z^z$ modulo reparametrizations (i.e. we erase the information where the root $z$ is located on $Z^z$).
Then, we can view the measure 
$M$ defined from $\hat M$ by 
$$ \frac {dM}{d\hat M} = \frac {1} {T(Z^0)}$$
 as a measure on the set of unrooted Brownian loops $Z$.
This is the Brownian loop measure.
%
%
%

Note that it is possible to conversely start with the measure $M$ on unrooted loops, and to recover
$\hat M$ by choosing a root for each loop $Z$, uniformly with respect to its time-parametrization.

%
It turns out that with this definition of $M$, if we define $M_D = M 1_{Z \subset D}$ to be the Brownian loop-measure restricted to those loops that stay in the domain $D$, then $M$ satisfies the conformal restriction property. More precisely:
\begin {lemma}[\cite {LW}]
\label {L4}
For any open domain $D \subset \C$, and for any conformal map $\Phi$ from $D$ onto $\Phi (D)$, 
$$ \Phi \circ M_D = M_{\Phi (D)}. $$
\end {lemma}
This fact is derived in \cite {LW}, see also \cite {Wln3}. It is a rather simple consequence of the 
conformal invariance of $N^z$. The idea is to start from the measure $M_D$, and to choose the root on a Brownian loop $Z$ with respect to another 
distribution than the uniform distribution. One chooses the root $z$ according to 
the uniform distribution with respect to the ``natural time'' of $\Phi(Z)$. Then, one just has to check that the 
couple $(\Phi(Z), \Phi (z))$  is indeed defined under $M_{\Phi(D)}$. 
\medbreak

Note that in this lemma, we did not require here the domain $D$ to be simply connected. However,  with our present definition, we restrict ourselves to the cases where $D \subset \C$ (i.e. it does not have a higher genus). 
We shall discuss the other cases later in the paper. 

\begin {figure}[t]
\begin {center}
\includegraphics [width=2.5in]{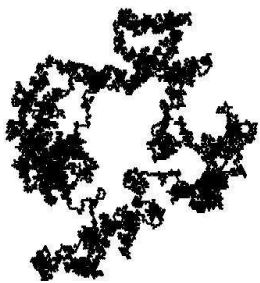}
\caption {The trace of a (discretized) Brownian loop (its outer boundary is Fig. {\ref {f1}})}
\end {center}
\label {f.3}
\end {figure}

\subsection {Outer boundaries}

Lemma \ref{L4} shows that
 $M$ does satisfy the same type of conformal restriction property as the one that we are looking for. But, of course,  Brownian loops are not self-avoiding. However, for each Brownian loop $Z$, we can define 
its ``outer boundary'' $\partial  Z$. More precisely, it is the boundary of the unbounded connected component of 
the complement of $Z$. It is not difficult to see that it is a continuous curve (it is the outer boundary of a continuous loop) and furthermore, it can be proved (e.g. \cite {BL}) that it is a self-avoiding loop (in other words, the Brownian loop has no cut-points --  these would correspond to double cut-points on a Brownian curve, and a priori 
estimates can be used to prove that such points do almost surely not exist on a Brownian curve).

Hence, $M$ induces a measure on (unrooted) self-avoiding loops in the plane. 
This measure exists and it satisfies conformal restriction in simply connected domains, because the outer boundary of the image of $Z$ under a conformal map 
$\Phi$ defined on a simply connected domain $D$, is the image under $\Phi$ of the outer boundary of $Z$. 
Combining this with the results of the previous section, we get the existence of the measure $\mu$ 
(note that this measure can be viewed as defined on ${\cal G}$ because $\hat M$ and $M$ are a priori defined on a richer $\sigma$-field; there are no measurability problems here): 

\begin {proposition}
There exists a non-trivial measure $\mu$ on the set of self-avoiding loops in the plane that does satisfy weak conformal restriction.
 It is equal (up to a multiplicative constant) to 
the measure on outer boundaries of Brownian loops defined via $M$.
\end {proposition}

Note that the conformal restriction property for $M$ works also for non-simply connected domains. 
However, this does not imply strong conformal restriction for the measure $\mu$ on self-avoiding loops; for instance, the image of the $\partial Z$ under 
the mapping $z \to 1/z$ is the ``inner boundary'' (i.e. the boundary of the connected component that contains the origin) of $1/Z$ and not its outer boundary. 
 
Let us now focus on the measure $\nu$ on loops surrounding the origin described in the previous section and 
on its Brownian interpretation. Clearly, taking
the measure $\mu$ restricted to those loops that surround the origin gives a way to construct $\nu$.

On the other hand, the measure $N^0$ on Brownian loops rooted at the origin defines also a measure on 
self-avoiding loops surrounding the origin. Furthermore, the conformal invariance property of $N^0$
implies that this measure on self-avoiding loops also satisfies the 
conformal restriction property that we require for $\nu$. Hence:

\begin {proposition}
\label {p6}
There exists a measure $\nu$ that satisfies the conditions of  Lemma \ref {L1}. It is equal (up to a multiplicative constant) to 
the measure on outer boundaries of Brownian loops defined under $N^0$. It is also equal (up to a multiplicative constant) to the measure on outer boundaries of Brownian loops surrounding $0$ defined under $M$. 
\end {proposition}

Let us insist on the fact that the two Brownian measures $M 1_{\{\partial Z \hbox{ surrounds } 0\} }$
and $N^0$ are very different, so that this lat result can seem somewhat surprising.
The second measure on loops is supported on rooted loops that go through the origin, 
while for the first one, (almost) no loop does go through the origin.
 
\subsection {Consequences for Brownian loops}

We now collect some further consequences of the previous facts concerning the Brownian loop measures.
These results are not really needed in our construction of the measure $\mu$ for general surfaces, but we 
write these consequences now to emphasize that their proof does (mostly) 
not use any SLE-based argument, but just properties of 
planar Brownian motion and conformal maps.

\medbreak

Proposition \ref {p6} has the following consequence:

\begin {corollary}
Consider two points $z$ and $\tilde z$ in the complex plane. 
The measure on outer boundaries of Brownian loops defined under $N^z$ and $N^{\tilde z}$ 
restricted to those loops that surround both $z$ and $\tilde z$ are identical.
\end {corollary}

This is just due to the fact that both measures coincide with (a constant times) 
the one defined by $M$ on this family of loops.  \qed

\medbreak

Let us remark that it is important in the previous corollary to consider the measures
$N^z$ with no prescribed time-length for the Brownian loops. Indeed, the following fact 
shows that some maybe tempting conjectures are false:

\begin {proposition}
Define a Brownian loop of time-length $T$ started from the point $z$, and its outer boundary $\gamma$.
Define another Brownian loop of time-length $T$ started from the point $\tilde z$ and its outer boundary $\tilde \gamma$. The laws of $\gamma$ and of $\tilde \gamma$ on the set of loops that surround both $z$ and $\tilde z$ do not (always) coincide.
\end {proposition}

\noindent 
{\bf Proof.}
One possible proof goes as follows:
Consider $z=0$, $\tilde z = x$ and $T = 1$, for very large real $x$. Note that loops of time-length $1$ that surround
both the origin and $x$ are rare. Then, we focus on the mass of the set of loops 
that surround the circle of radius $x$ around the origin. 
In the case of the Brownian loop from $x$ to $x$, a typical loop of time-length $1$ with an outer boundary that 
surrounds the unit circle is going very fast around the circle. In the case of the loop from $0$ to $0$, 
the loop has first to reach the circle, then to go around it, and finally to come back to the origin, and this 
has a much smaller probability (in the $x \to \infty$ regime).  
\qed
   
\medbreak

Let us make a little general comment
 on shapes for scale-invariant measures on random compact sets. We are going to 
write it in the case of $N^0$, but it can clearly be adapted to the other measures that we will discuss in this 
paper. 
If we ``scale out'' the scale factor, we can see that $N^0$ defines a probability measure on the set of possible ``shapes''. In this Brownian case, we say that $Z$ and $Z'$ have the same shape if for some 
$\lambda>0$, $Z= \lambda Z'$. In other words, we look at the set of rooted loops modulo scaling.

There exists many different (measurable) 
ways to choose a given representative for each shape (i.e. for each equivalence class of loops). 
One can for instance 
use the ``natural Brownian time-parametrization'' of the loop $Z$, and choose the representative such that $T=1$
(this definition works well for $N^0$-almost all loops).
Another possibility is to choose the representative such that the maximal distance to the origin on the loop is $1$.

Consider the probability measure $P_{0\to 0; 1}$ on Brownian loops in the plane of time-length $1$ (i.e. the real
and imaginary parts are independent Brownian bridges of time-length $1$). Define the 
product measure, $dr/r \otimes P_{0 \to 0; 1}$. If we look at $rZ$ where $(r,Z)$ is defined under this product measure, it is easy to see that it is defined under the measure $N^0$.
Similarly, scaling shows that 
$$N^0 = \int_0^\infty \frac {dT}{2T} P_{0 \to 0 ; T },
$$
where $P_{0 \to 0 ; T}$ is this time the law of a Brownian loop of time-length $T$.

Clearly, for any positive real random variable $U$ (that can depend on $Z$ -- but does not interact with $r$), we see that the measure of $rU Z$ is always the same. In particular, this shows that one can view 
$N^0$ as obtained by taking the product of $r$ (defined under the measure $dr/r$) and a probability measure
$P$ on the space of representatives (i.e. of shapes) of Brownian loops. 

Conversely, one way to define $P$ from $N^0$ goes as follows (we safely leave this as an exercise to the reader):
 Define the set of loops $L_r$ that lie in the annulus 
$\{z \ : \  1< |z| < r \}$. Then, consider the measure on shapes defined by 
$N^0( \cdot  \cap L_r) / N^0 (L_r)$. This is clearly a probability measure. Consider then the 
limit of this probability measure as $r \to \infty$. Note that $N^0$ and $2N^0$ give rise to the same probability 
measure on shapes.

Hence, in a way, defining the scale-invariant measure $N^0$ (modulo a multiplicative constant) is equivalent to defining the probability measure $P$ on shapes of loops.
 
\medbreak
  
Looking at Proposition \ref {p6},
it is of course natural to ask what the multiplicative factor relating the outer boundaries defined under the
 measures $M$ and $N^0$  is, and to what constant $c$'s these measures correspond to in the formula (\ref {E4}). 
In the case of $N^0$, things are easy: We have chosen our normalization in such a way that the
$N^0$ mass of the set of loops with maximal radius between $1$ and $e$ is $1$. Hence, it follows immediately that the constant $c$ in (\ref{E4}) is equal to 1 for $N^0$.

In the case of $M$, it is more complicated to compute the constant $c$.
The starting point of the loop 
is chosen away from the origin and one has to keep only those loops that surround the origin.
Let us now look at the set $A$ of self-avoiding loops that surround the origin and have time-length 
between $1$ and $e^{2r}$.
Clearly,
\begin {eqnarray*}
M(A) &=&
\int_\C d^2 z N^z ( \frac {1_A}{T(Z)}  ) \\
&=& \int_\C d^2 z \int_1^{\exp (2r)} \frac {dT}{2T^2} P_{z \to z ; T} ( Z \hbox { surrounds } 0 ) \\
& = &\int_\C d^2 z  \int_1^{\exp (2r)} \frac {dT}{2T^2} P_{0 \to 0 ; T} ( Z \hbox { surrounds } z) \\
&=& \int_1^{\exp (2r)} \frac {dT}{2T^2} E_{0\to 0; T} ( {\cal A} (Z)) dT ,
\end {eqnarray*}
where ${\cal A}(Z)$ denotes the area of the inside of the self-avoiding loop defined by $Z$. 
Clearly, 
$ E_{z\to z; T} ( {\cal A} (Z))
=T E_{z\to z; 1} ( {\cal A} (Z))$ because of scaling.
Furthermore, we know from \cite {GT} that 
$$E_{z\to z; 1} ( {\cal A} (Z)) = \pi /5 .$$
Note that this result has been obtained in \cite {GT} 
via considerations involving conformal restriction and SLE (this is the only time in this section where we use an SLE-based result).
Hence,
$$ M(A) = \frac {\pi}{5} r. $$
Note finally that the measure $M$ restricted to those Brownian loops that surround the 
origin is scale-invariant. Hence, at each scale, the area ${\cal A}(Z)$ scales like the time-length of $Z$.
It therefore follows readily that $c$ for the case of the measure on self-avoiding loops surrounding the 
origin induced by $M$ is equal to $\pi /5$, and therefore this measure is equal to $\pi /5$ times the one defined via $N^0$. 

\medbreak

In fact, it is possible to construct various other multiples of $\nu$ using the Brownian loop measure $M$. We list some of the possibilities:
\begin {itemize}
\item
Note that a Brownian loop $Z$ defines an integer index around any point $z \notin Z$, and in particular, around the origin (since for $M$-almost every loop $Z$, $0 \notin Z$). Let us call ${\cal N}$ this index.
Then, for each fixed $n \in \Z$, we define the measure on self-avoiding loops surrounding the origin induced by $M$, but restricted to those loops $Z$ such that ${\cal N} = n$.
It is easy to see that this measure satisfies conformal restriction too, so that it is a multiple of the measure $N^0$.
The constant can be easily identified via the law of the index of the Brownian loop (see \cite {Y,GT}).
\item
For each loop $Z$, we define the connected components $C_j$ of the complement of $Z$. We define a graph on this set of 
connected components as follows: We say that $C_j$ and $C_i$ are neighbors, if there exists a point on the loop 
that touches both $\partial C_i$ and $\partial C_j$. It is easy to see that the existence of cut-points on planar Brownian curves (first proved in \cite {Bu}) implies that each $C_j$ has infinitely many neighbors.
But it is an open question whether this graph has several connected components or not.   
Let us call $C_0$ the connected component containing the origin, and $C_\infty$ the connected component containing infinity. For each $n \ge 1$, we can define the event $A_n$ that $C_0$ and $C_\infty$ are at distance exactly $n$ in the graph. If we now focus on the outer boundaries of Brownian loops restricted on the set of loops for which $A_n$
holds, we get also a measure satisfying conformal restriction. Hence it is also a multiple of $N^0$.
This time, it seems difficult to compute the multiplicative constant.
\end {itemize}

At this point, it is probably worth mentioning that the Brownian loop measure rooted at $z$ and
the unrooted Brownian loop measure defined in \cite {LW} and denoted there by $\mu(z,z)$ and $\mu^{loop}$ differ from
$N^z$ and $M$ by a multiplicative constant. More precisely, $\mu(z,z) = N^z / 2\pi$ and 
$\mu^{loop} = M / 2\pi$. The normalization there is due to the relation between these measures
and the ``central charge'' of the corresponding models, in the convention used in theoretical physics. In the normalization of \cite {LW}, a Poissonian cloud of intensity $c \mu^{loop}$ (see \cite {Wls}) correspond to a representation with central charge $c$.

\medbreak

To conclude this section,
 we do not resist the temptation of citing the following visionary phrases from Mandelbrot \cite {Ma}:

\medbreak
``{\sl a Brown loop separates the plane into two parts: an exterior which can be reached from a distant point without 
intersecting the loop, and an interior which I propose to call a Brown hull or Brown Island. (...) I propose for the 
Brown hull's boundary the term self-avoiding Brownian motion.}''
  
\section {Annular regions}

We have now constructed a measure on self-avoiding loops  as outer boundaries of 
planar Brownian loops. 
This measure satisfies weak conformal restriction and the construction of the measure is 
very non-symmetric: The ``outside'' boundary of the Brownian loop is (by definition) 
on the same side of the Brownian path
 as the boundary of the simply connected domain that 
 we are looking at. It seems that there is little hope to 
say anything similar in annular regions, since the notion of inside and outside would become symmetric.

One of the main points of the present paper is however to show that -- even if its Brownian construction does not work anymore -- the measure on self-avoiding loops that it defines does still satisfy conformal restriction in annular regions, and therefore also on any Riemann surface.

The main statement of the present section is the following:

\begin {proposition} 
\label {P4}
Consider two annular regions $D$ and $\tilde D$ in the complex plane, such that there exists a conformal map 
$\Phi$ from $D$ onto $\tilde D$ (i.e. the two regions have the same modulus). 
Define $\mu_D$ (respectively $\mu_{\tilde D}$) as the measure $\mu$ restricted to those loops that stay in $D$ (resp. 
$\tilde D$). Then, $\Phi \circ \mu_D = \mu_{\tilde D}$.
\end {proposition}

Let us insist on the fact that in general, such a conformal map $\Phi$ does not extend conformally to the inside nor to the outside of the annular region, so that this is a much stronger conformal invariance statement than weak conformal restriction.

Note that the fact that these two measures coincide on the set of self-avoiding loops that do not go ``around'' the hole is not surprising and follows readily from the result in the simply connected case. However, for those loops that 
go around the hole, the result can seem surprising, and as we shall see later, it has some non-trivial 
consequences for Brownian loops.

In order to prove this proposition, we first need to collect and recall some relevant facts from \cite {LSWr, LW, Wcrr}.

In \cite {LW}, we showed how to decompose a Poissonian cloud of Brownian loops in a simply connected
domain, sampled according to the intensity $M$, according to the first point on a Loewner chain that the 
Brownian loops intersect. We are now going to recall this result, but in its ``non-Poissonized'' version 
(looking at the intensity measure instead of the Poissonian sample).

Consider the upper half-plane $\H$. Suppose that $\eta: [0,T] \to {\overline \H}$ is a continuous curve 
without double-points in $\H$, such that $\eta (0)= 0$, $\eta (0,T] \in \H$. 
We are going to suppose that $\eta$ is parametrized in such a way that for each $t$, the
conformal map $\varphi_t$ from $\H \setminus \eta(0,t]$ onto $\H$ such that $\varphi_t (\eta_t)= 0$,
$\eta_t(z) \sim z$ as $z \to \infty$ satisfies 
$$ \varphi_t (z) = z + b + \frac {2t}z + o (1/z) 
$$
as $z \to \infty$ for some real constant $b=b(t)$. For any simple continuous 
path $\eta$, such a parametrization clearly exists
(in other words, $\eta$ is a Loewner chain, parametrized by its capacity from infinity).

We need to define the Brownian bubble measure $Bub$ in the upper half plane:
This is a measure supported on Brownian loops starting and ending at the origin, but that stay in $\H$ 
in the meantime. Very roughly speaking, one can view $Bub$ as the ``a version of $N^0$ conditioned on the 
loops that stay in $\H$'' (or a version of $M$ conditioned on the set of loops that touch the boundary of $\H$ just at $0$ and stay in $\H$ otherwise).
More precisely,
let $Q_\eps$ denote the law of a Brownian path started from $i \eps$ and conditioned to exit the 
upper half-plane in the interval $[0,\eps]$ (and stopped at this exit time).
Note that the $Q_\eps$ probability that the Brownian path reaches the circle of radius $1$ is 
of the order $\eps^2$. Indeed, the conditioning in the definition of $Q_\eps$ is with respect to an event of constant probability (independent of $\eps$) because of scaling. On the other hand, a Brownian motion started from 
$i\eps$ reaches distance one from the origin without exiting the positive half-plane with probability of the 
order $\eps$. Then, it has also a conditional probability $\eps$ to exit $\H$ via the small interval $[0, \eps]$
of size $\eps$. Hence, it should not be really surprising that
$$ Bub = \lim_{\eps \to 0} \eps^{-2} Q_\eps $$
exists and that it has the following scaling property: If $Z$ is defined according to $Bub$, then 
$\lambda Z$ is defined according to $\lambda^{-2} Bub$.
For more details, see \cite {LSWr,LW}.

 The following lemma is 
derived in \cite {LW} in its poissonized form.
\begin {lemma}[\cite {LW}]
\label {L2}
For some positive constant $c$,
$$
M_\H ( \cdot 1_{ Z \cap \eta [0,T] \not= \emptyset } )
=
c \int_0^T dt \varphi_t^{-1} ( Bub (\cdot ))
.$$
\end {lemma}

This lemma should not be surprising: We are decomposing the measure $M_\H$ restricted to those Brownian loops 
that intersect $\eta$, according to the smallest $t$ such that $\eta_t \in Z$. And, intuitively speaking:
\begin {itemize}
\item
Because of conformal restriction, the conformal Markov-type property (that lead to the definition of SLE in 
a slightly different context) holds, so that 
$$M_\H ( \cdot 1_{ Z \cap \eta [0,T] \not= \emptyset } )
=
 \int_0^T dt \varphi_t^{-1} ( Bib (\cdot ))
$$
for some measure $Bib$ supported on the set of bubbles in $\H$.
\item
This measure $Bib$ is equal to a constant times $Bub$, because when $T$ is very small, it is very close to the 
``conditioned'' version of $M$.
\end {itemize}
We do not care here about the values of the constants as they are not needed for our purpose.

\medbreak

We now derive another important intermediate fact, and this is where considerations involving SLE$_{8/3}$ come into play in 
a crucial way, even if it is not so apparent. The SLE-ingredients that are used are that 
SLE$_{8/3}$ is a simple curve (this is proved in \cite {RS})  that
 satisfies two-sided conformal restriction in the chordal case (this is proved in \cite {LSWr} and briefly recalled at the end of the introduction) :

Let us focus on the outer boundary $\gamma$ of a Brownian loop
defined under the Brownian bubble measure $Bub$. We denote by $bub$ the measure under which $\gamma$
 is defined. 
It turns out (e.g. \cite {LSWr}) that one define this measure using SLE:

\begin {lemma}[\cite {LSWr}]
\label {l11}
The measure $bub$ is equal to
the limit when $\eps \to 0$ of 
$$ cst \times \eps^{-2} \times P_{0 \to \eps}$$
where $P_{0 \to \eps}$ denotes the law of the chordal SLE$_{8/3}$ curve from $0$ to 
$\eps$ in the upper half-plane.
\end {lemma}

 In the sequel, 
$P_{x \to y}$ will denote the law of chordal SLE from $x$ to $y$ in $\H$.  
Consider now a simply connected bounded closed set $A \subset \H$, and a conformal map $\Phi$ from the annular region 
$\H \setminus A$ onto some other annular region $\H \setminus \tilde A$ such that $\Phi (\infty) = \infty$.
For convenience, we define 
$$ k = k (\eps, \Phi) = 
\left(
\frac {\eps^2 \Phi'(0)\Phi'(\eps)}{(\Phi(\eps) - \Phi(0))^2} 
\right)^{5/8}.$$

\begin {lemma}
\label {L10}
Consider an SLE$_{8/3}$ curve $\beta$ from $0$ to $\eps$ in $\H$. The law of 
$\Phi(\beta)$, restricted to the event where $\beta \subset \H \setminus A$ 
is identical to $k_\eps$ times the law of an SLE$_{8/3}$ from 
$\Phi(0)$ to $\Phi(\eps)$ in $\H$, restricted to stay in $\H \setminus A$.
\end {lemma}
In other words,
$$
\Phi ( 1_{\beta \subset \H \setminus A} P_{0 \to \eps} ) 
= k (\eps, \Phi) P_{\Phi(0) \to \Phi(\eps)} 1_{\beta \subset \H \setminus \tilde A}.
$$
A related  fact has been derived by Vincent Beffara in his PhD thesis \cite {Bephd}, section 5.1.2.
We give here a direct justification based on SLE$_{8/3}$'s chordal conformal restriction property:

\medbreak
\noindent
{\bf Proof.}
Consider a simple path $\delta$ that joins the real line to $A$ in $\H$ (and does start away from $0$ and $\eps$
on the real line). Note that the starting point of $\delta$ can be in the interval $(0, \eps)$.
Consider the image $\tilde \delta$ of $\delta$ under $\Phi$.  Define the simply connected sets 
$H=\H \setminus (\delta \cup A)$ and $\tilde H= \H \setminus (\tilde \delta \cup \tilde A)$. 
$\Phi$ maps conformally the former to the latter. Hence, by conformal invariance, 
the image of SLE$_{8/3}$ from $0$ to $\eps$ in $H$ is an SLE$_{8/3}$ from 
$\Phi(0)$ to $\Phi(\eps)$ in $\tilde H$.
We would like to say that the law of $1_{\beta \subset H} \Phi \circ \beta$ is identical to 
$ k P_{\Phi(0) \to \Phi(\eps)} 1_{\beta \subset \tilde H}$. Since by chordal restriction 
and conformal invariance,  both these measures are in fact
multiples of the SLE law in $\tilde H$, it remains only to compare the total masses.

The probability that an SLE$_{8/3}$ from $0$ to $\eps$ in $\H$ stays in $H$ is $(\psi'(0) \psi'(\eps))^{5/8}$, 
where $\psi$ denotes a conformal 
map from $H$ onto $\H$ that fixes $0$ and $\eps$. This is SLE$_{8/3}$'s chordal restriction property \cite {LSWr}.

Similarly, the probability that an SLE$_{8/3}$ from $\Phi(0)$ to $\Phi(\eps)$ in $\H$ stays in $\tilde H$, is
$(\tilde \psi' (\Phi(0)) \tilde \psi' (\Phi (\eps)))^{5/8}$, where  $\tilde \psi$ is a conformal map from $\tilde H$ onto $\H$ that fixes $\Phi(0)$ and $\Phi(\eps)$.

Note that
$ \psi^{-1} \circ \Phi \circ \psi$ maps conformally  $\H$ onto itself, and that it maps $0$ to $\Phi (0)$ 
and $\eps$ to $\Phi(\eps)$. Hence, the product of its derivatives at $0$ and at $\eps$ is equal to
$(\Phi(\eps)-\Phi(0))^2 / \eps^2$ (because a simple scaling does the job). 

Putting the pieces together, 
we get that for any $\delta$, the law of 
$$1_{\beta \subset \H \setminus (A \cup \delta)}  \Phi \circ \beta$$ 
is indeed identical to 
$k P_{\Phi(0) \to \Phi(\eps)} 1_{\beta \subset \H \setminus (\tilde A \cup \tilde \delta)}$.
Since this is true for any $\delta$ (and because a simple curve from $0$ to $\eps$ in $H$ does not disconnect $A$ from the real axis), the lemma follows. \qed

\medbreak

We now write down the consequence of the above for the bubble measure $bub$ on
self-avoiding loops in the upper half-plane:

\begin {lemma} 
Let $A$, $\tilde A$ and $\Phi$ be as in Lemma \ref {L10}. We also assume that $\Phi(0)= 0$.
Then, the measure $\Phi \circ (bub  1_{\gamma \cap A = \emptyset})$
is identical to 
$$ \Phi'(0)^{-2} bub 1_{\gamma \cap \tilde A  = \emptyset}.$$
\end {lemma}

\noindent
{\bf Proof.}
This follows readily from Lemma \ref {l11}  and the properties of the SLE measure described in the last lemma. Note also of course that 
$$ \lim_{\eps \to 0} \frac {\eps^2 \Phi'(0)\Phi'(\eps)}{(\Phi(\eps) - \Phi(0))^2} = 1. $$ 
\qed

\medbreak

We are now ready to prove the proposition:

\medbreak
\noindent
{\bf Proof of Proposition \ref {P4}.}
Consider first two annular regions of the type $\H \setminus A$ and $\H \setminus \tilde A$
such that there exists a conformal map $\Phi$ from the first onto the second.
It suffices to prove that for any Loewner chain $\eta [0,T]$ in $\H$ that does not intersect $A$, 
the image under $\Phi$ of the measure $\mu$ restricted to those loops in $\H$ that do intersect $\eta$
but not $A$, is exactly identical to the measure $\mu$ restricted to those loops in $\H$ that do 
intersect $\Phi(\eta)$ but not $\tilde A$.

\begin{figure}[t]
\begin{center}
\includegraphics[height=2in]{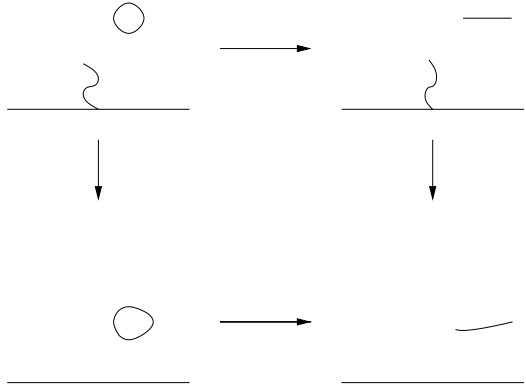}
\caption{The commutative diagram with $\Phi$, $\Phi_t$, $\varphi_t$ and
  $\tilde \varphi_s$}
\label{figure:example}
\end{center}
\end{figure}

We know that 
\begin {equation}
\label {eq14}
1_{\gamma \cap A = \emptyset, \gamma \cap \eta [0,T] \not= \emptyset} \mu_\H 
=
\int_0^T  dt \varphi_t^{-1} ( bub  1_{\gamma \cap \varphi_t (A) = \emptyset} ) 
.\end {equation}
Hence, the conformal image under $\Phi$ of the left-hand side of the previous expression is equal to 
$$
\int_0^T  dt \Phi \circ \varphi_t^{-1} ( bub ( \cdot) 1_{\gamma \cap \varphi_t (A) = \emptyset} )
.$$
The path $\Phi (\eta [0,T])$ is also a Loewner chain. But its parametrization by capacity is different from
that of $\eta$. Let $s = s(t)$ denote the capacity at infinity of $\Phi ( \eta [0,t])$ in $\H$, and 
let us define $\tilde \eta$ by $\tilde \eta_{s(t)} = \Phi ( \eta_t)$, so that $\tilde \eta$ is a 
naturally parametrized Loewner chain. 
Let us define the conformal map
$\tilde \varphi_s$ from $\H \setminus \tilde \eta [0,s]$ onto $\H$, such that $\tilde \varphi_s (z) 
\sim z$ at infinity, and $\tilde \varphi_s (\tilde \eta_s) = 0$.

We now define $A_t = \varphi_t (A)$, $\tilde A_s = \tilde \varphi_s (\tilde A)$, and $\Phi_t$, the conformal map from
$\H \setminus A_t$ onto $\H \setminus \tilde A_{s(t)}$ defined by 
$$
\Phi_t = \tilde \varphi_{s(t)} \circ \Phi \circ \varphi_t^{-1}
.$$
It is normalized at infinity, and maps the origin to the origin.
Note that 
$$
 s(t) = \int_0^t du |\Phi_u' (0)|^2 
$$ 
for scaling reasons: The capacity of $\varphi_t (\eta[t, t + \eps])$ is $\eps$ by definition. The capacity of 
$\tilde \varphi_{s(t)} ( \tilde \eta [ s(t) , s(t + \eps) ])$ is $s(t+\eps) - s(t)$ by definition. But 
the second slit is the image of the first one by $\Phi_t$, and it is then standard that 
$s(t+  \eps) - s(t) \sim \eps |\Phi_t'(\eta_t)|^2$ when $\eps \to 0$.

Hence, if $S=s(T)$,
\begin {eqnarray*}
\lefteqn { 
\Phi \circ (1_{\gamma \cap A = \emptyset, \gamma \cap \eta [0,T] \not= \emptyset} \mu_\H ) }\\
&=&
\int_0^S  ds |\Phi_t' (0)|^{-2} 
 \tilde \varphi_s^{-1} \circ \Phi_t ( bub ( \cdot) 1_{\gamma \cap \varphi_t (A) = \emptyset} )\\ 
&=&
\int_0^S  ds  
 \tilde \varphi_s^{-1} \circ bub ( \cdot) 1_{\gamma \cap \tilde A_s = \emptyset} )
\\
&=&
1_{ \gamma \cap \tilde A = \emptyset, \gamma \cap \tilde \eta [0,S] \not= \emptyset }
\mu_{\H}.
\end {eqnarray*}
But this is true for any Loewner chain $\eta[0,T]$ in $\H \setminus A$. Hence, it follows that
$$ \Phi ( 1_{\gamma \cap A= \emptyset} \mu_\H ) = 1_{\gamma \cap \tilde A = \emptyset} \mu_\H.$$
This is 
exactly the statement of the Proposition in the case where $D= \H \setminus A$ and $\tilde D= \H \setminus \tilde A$.

In the general case (i.e. general annular domains $D$ and $\tilde D$), define $D_1$ (resp. $\tilde D_1$) the complement of the unbounded connected component of $\C \setminus D$  (resp. $\C \setminus \tilde D$) i.e. the domain obtained by filling the hole of $D$ (resp. $\tilde D$). 
One can map conformally $D_1$ (resp. $\tilde D_1$) onto the upper half-plane, and apply the result to the images of 
$D$ (resp. $\tilde D$) under these maps (these are sets of the type $\H \setminus A$ and $\H \setminus \tilde A$ as before). Then, mapping things back onto $D$ and $\tilde D$, we get immediately the proposition for $D$ and $\tilde D$.  \qed

\section {General Riemann surfaces}

We can now turn our attention to general Riemann surfaces. In the sequel, we suppose that $\mu$ is defined in such a way that the constant $c$ defined in (\ref {E4}) is equal to one.

\subsection {Definition of $\mu$ and conformal restriction}

 Let $\rho (D)$ denote the modulus of an annular region $D$  i.e. the number such that 
there exists a conformal map from $D$ onto the annulus $\{ z \ : \ 1 \le |z| < \exp (\rho )\}$. 
We say that a self-avoiding loop in $D$ goes around the hole in $D$, if it disconnects the inner part
of $\partial D$ from its outer part (in $D$).
Then: 

\begin {lemma}
There exists a function $F: (0, \infty) \to (0, \infty)$ such that for any annular region $D \subset \C$, 
$$ 
\mu ( \{ \gamma \ : \ \gamma \subset D \hbox { and goes around the hole in }D \} ) 
= F( \rho (D)).
$$ 
\end {lemma}

This is just due to the fact that $\mu$ satisfies conformal  invariance for annular regions. \qed

\medbreak

This will allow us to justify the following definition of $\mu_S$ for any Riemann surface $S$:
 
\medbreak
\noindent
{\bf Definition.}
{\sl Let us now consider any Riemann surface $S$. We say that $D \subset S$ is an annular region if it is conformally equivalent to an annulus. Then, we { define} $\mu_S$ by the fact that for any annular region $D \subset S$,
$$ 
\mu_S ( \{ \gamma \ : \ \gamma \subset D \hbox { and goes around the hole in } D \} ) 
= F( \rho (D)),
$$ 
where $F$ is the function defined in the previous lemma.}
\medbreak

Let us first show that this definition indeed determines $\mu_S$ uniquely:
The family of events of the type
$$A_D=  \{ \gamma \ : \ \gamma \subset D \hbox { and goes around the hole in }D \} $$
when $D$ varies in the family of annular regions in $S$
 is stable under finite intersections.
It is also easy to check that this family of events generates the 
$\sigma$-field on which $\mu_S$ is defined.
Hence, standard measure-theoric considerations imply that there exists at most one measure 
$\mu_S$ such that $\mu_S (A_D) = F ( \rho (D))$ for any annular region $D \subset S$.

Let us now check that the measure $\mu_S$ exists.
Each self-avoiding loop in $S$ is contained in some annular domain $D \subset S$ (consider
for instance  a ``small neighborhood'' of the loop). One can in fact easily define a countable 
family of annular domains $D_n \subset S$ such that each loop is contained in at least one domain $D_n$. 
We have defined for each $D_n$ a measure $\mu_{D_n}$ on loops in the annular region $D_n$ in the previous 
section. These families are compatible since 
$\mu_{D_n}$ and $\mu_{D_m}$ do coincide on $D_n \cap D_m$ (the measure do satisfy restriction 
in annular regions). Hence, this indeed constructs a measure $\mu_S$ on loops in $S$, and  for any $n$, $\mu_S  1_{\gamma \subset D_n} = \mu_{D_n}$.

The conformal restriction property of the family of measures $\mu_S$ then follows immediately from the conformal restriction property in the annular regions derived in the previous section.
\medbreak

We have now in fact derived all the ingredients needed to check Theorem \ref {T1}:

\medbreak
\noindent
{\bf Conclusion of the proof of Theorem \ref {T1}.}
We have just constructed a family $(\mu_S)$ that satisfies conformal restriction.
It remains to check that it is the unique one (up to a multiplicative constant).

To prove this, we can note that the previous argument shows that the
whole family $(\mu_S)$ is characterized by the subfamily $(\mu_D)$ where $D$ spans the family of 
annular regions. Hence, the whole family is in fact characterized by $\mu=\mu_\C$ because one can then
define all $\mu_D$'s for annular regions by restriction.
But we have seen that there exists a unique measure $\mu$ (up to a multiplicative constant) on loops in 
$\C$ that satisfies conformal restriction for simply connected domains. Hence, the family $(\mu_S)$ is unique too.
\qed

\subsection {Relation to Brownian loops, zero genus.}

Let us recall the steps in our construction of $\mu$. In the simply connected case (and for genus $g=0$), we constructed 
the measure $\mu$ using outer boundaries of Brownian loops defined under the loop measure. Then, using the relation 
with SLE$_{8/3}$, we proved that this measure $\mu$ satisfies an inside/outside symmetry property, and this enabled 
us to define the measure $\mu$ on any Riemann surface. This raises naturally the question of the relation between 
the measure $\mu$ on a surface $S$ and measures on Brownian loops on $S$.

Let us first focus on the case where the genus of $S$ is zero, i.e. $S=D$ is a subset of the Riemann sphere.
In this case, the Brownian loop measure in $D$ is simply $M_D$ i.e. the restriction of the measure $M$ on Brownian loops
in the plane to the set of loops that stay in $D$.

If one considers a bounded annular region $D \subset \C$ such that the origin 
is in the hole, we can consider the conformal map $z \mapsto 1/z$ defined from 
$D$ onto the annular region $1/D$. The image of $\mu_D$ is exactly $\mu_{1/D}$.
Hence, the image of the self-avoiding curve $1/\gamma$ can be viewed 
as the outside boundary of a Brownian loop $Z'$. But because of the invariance of the Brownian loop-measure $M$
under the map $z \mapsto 1/z$ (this is due to its conformal invariance), 
we see that we can also view $1/\gamma$ as the ``inside boundary'' of a Brownian loop $1/Z$.
In fact, we have just proved the following:

\begin {proposition}
 Consider a Brownian loop $Z$ defined under the Brownian loop measure $M$, 
 and call its outer boundary 
$\gamma$ as before. Define also the boundary $\bar \gamma$ of the connected component of $\C \setminus Z$ that 
contains the origin. Note that $\gamma = \bar \gamma$ unless $Z$ disconnect $0$ from infinity.
Then, the measure under which $\bar \gamma$ is defined is identical to the measure under which 
$\gamma$ is defined (and is equal to $\mu$).
\end {proposition}

\medbreak

Suppose that $Z$ is a Brownian loop defined under the measure $M$. Define the family $(\gamma_j)$ of its boundaries (i.e. the boundaries of the connected components of its complement. 
Then, define the measure on self-avoiding loops $\hat\mu$ as 
$$ \hat \mu  ( A) = M ( \sum_j 1_{ \gamma_j \in A } )$$ 
(i.e. we take all loops simultaneously). Note that the complement of the Brownian loop has of course only finitely many large connected components 
(the number of small loops has been studied in \cite {Mo, LG}). 

\begin {figure}
\begin {center}
\includegraphics[width=2.5in]{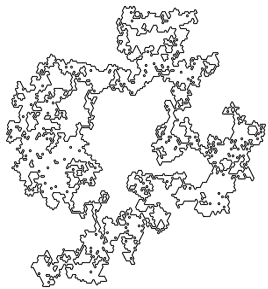}
\caption {The boundaries in the Brownian loop from Fig. 
3}
\end {center}
\end {figure}

\begin {corollary}
\label {c15}
We have $\hat \mu = 2 \mu$.
\end {corollary}

\noindent
{\bf Proof.}
The measure $\hat \mu$ is clearly translation-invariant. Furthermore, if we look at the measure $\hat \mu$ restricted to those loops that disconnect the origin from infinity, we see that we count once the inner boundary, and once the outer boundary (if $Z$ disconnects the origin) and nothing else. The result follows.
\qed

\medbreak
For any multiply connected subset $D$ of the plane, it is always possible to 
view $\mu_D$ as a measure on outer boundaries of Brownian loops: One just has to fill in all
the holes except a given prescribed one (this defines a simply connected subset of the Riemann sphere) and to consider outer boundaries of 
Brownian loops defined in this domain.

\medbreak

The conformal invariance of $\mu_D$ has some other consequences. For instance, 
consider 
$D_1 = \C \setminus \U$ and $D_2 = \C \setminus [-1,1]$. These two sets are conformally equivalent (let $\psi$ denote the conformal map from $D_2$ onto $D_1$ that is symmetric with respect to the real axis). 
One can view $\mu_{D_1}$ as the measure on outer boundaries $\gamma$ 
of Brownian loops $Z$ in $\C$, restricted to those loops such that $\gamma \cap \U = \emptyset$. 
But, one can also view it as the conformal image via $\psi$ of the measure of outer boundaries of Brownian
loops $Z$ in $\C$ such that $\gamma \cap [-1,1] = \emptyset$. 
Hence, as for outer boundaries that do not intersect $\U$ are concerned, one can replace the Brownian loops 
in $\C$ by ``Brownian loops in the simply connected domain obtained from $\C \setminus \U$ by identifying all conjugated points  $e^{i \theta}$ and $e^{- i \theta}$ on $\partial \U$.''

\medbreak
Let us also note that the fact that one Brownian loop defines infinitely many boundaries is not in contradiction 
with the fact that $\hat \mu = 2 \mu$ because we are talking about scale-invariant (infinite) measure. In fact, since 
the Lebesgue measure of a Brownian path is zero, it is clear that the area enclosed by the outer boundary of 
a Brownian loop is identical to the sum of all areas of the inner boundaries of the same loop (so that this sum converges), which is another way 
to see why this identity holds (with the factor $2$).

Also, recall that for a Brownian loop of time-length $1$, it is known \cite {Mo,LG} that the number of 
inner boundaries that enclose an area at least $u$ behaves almost surely like 
$2 \pi / ( u (\log (1/u)^2))$ as $u \to \infty$ (this question was raised by ... Mandelbrot).
Furthermore, the result of \cite {Wccc} shows that in fact it is enough to look at the sample of
one Brownian loop to describe fully the measure $\mu$ in $\C$: If $F$ is a function of the ``shape'' (i.e. a function invariant under scaling and translations) of a connected component, and if $C_1, C_2, \ldots$ are the bounded connected components of the complement of the Brownian loop
ordered by decreasing area, then  \cite {Wccc}
$ (F(C_1) + \ldots + F(C_n))/n $ tends (in $L^2$) to ``$P(F)$'', where 
$P (F)$ is in fact the expected value of $F$ for the probability measure on the 
shape of the outer boundary of a Brownian loop.

\subsection {Relation to Brownian loops: Higher genus}

Let us now say a few words on the relation between Brownian loop measures and $\mu_S$, when the surface $S$ has a higher genus (for instance when $S$ is the torus $\C / ( \Z + \tau \Z)$).

We can start with the case where $S$ is compact (since the other cases will follow from restriction).
We first have to define properly the Brownian loop measure. It is quite clear how to proceed:
For each $z \in S$, we can define the Brownian loop measure $N_S^z$ in $S$ rooted at $z$
in a similar way as in $\C$ by taking the limit when $\eps \to 0$ of $4\log (1/\eps)$
times the law of a Brownian motion on $S$, started uniformly on the circle of radius $\eps$ around  
$z$ and stopped at its first hitting time of the circle of radius $\eps^2$ around $z$. 
Note that we can define for each loop, its natural time $T(Z)$ (given for by its total quadratic variation), but that this notion depends on the choice of a metric on $S$ (they are only defined locally modulo conformal invariance).

After arguing that the map $z \to N^z$ is in some sense measurable (one can for instance couple 
$N^z$ and $N^{z'}$ when $z$ and $z'$ are close), we can define the measure
$$
\hat M_S = \int_S d^2 z N^z \frac {1}{T(Z)}
$$
where $d^2 z$ denote here the chosen metric on $S$. This is now a measure on the set of rooted Brownian loops in $S$. Finally, we erase the information of where the 
root is on $Z$, and define the corresponding measure $M_S$ on unrooted Brownian loops in $S$. Note that when the genus of $S$ is zero, this is the same as the Brownian measure that we have studied earlier in this paper. 
The proof of conformal invariance and restriction for this family of measures $M_S$ is then essentially identical to that of the restrictions of the planar measure $M_\C$. 

A first remark is that the boundaries of connected components of the complement of a Brownian loop defined under $M_S$ are not necessarily self-avoiding anymore, even for ``planar connected components''. Consider for instance the torus $S= \C / ( \Z + i \Z)$. It can happen (i.e. with positive mass) that the 
Brownian loop goes around the torus in both directions (i.e. for instance that in the covering surface $\C$ it goes through the neighborhood of the origin, then stays close to the segment $[0,1]$ and then to the segment $[1, 1+i]$ and closes the loop near $1+i = 0$) and that the path has local cut-points.
At such a cut point, locally, both sides of the Brownian motion are in the same connected component of 
$S \setminus Z$. Hence, the boundary of this connected component has a double point (and in fact, it has many cut points i.e. the dimension of the set of cut points is $3/4$, see \cite {LSW2}).

On the other hand, for a given surface $S$, the $M_S$-mass of the set of loops such that the boundaries of the connected components of its complement in $S$ are not all simple curves is finite. Hence, focusing our 
attention on those loops that have simple boundaries is not so restrictive.
One could then see whether there is  a  direct relation to the measure $\mu$.

The answer does not seem a priori obvious. Recall for instance from \cite {WToulouse} that a Brownian excursion 
``conditioned to have no cut-point'' defines a restriction measure of parameter $2$, i.e. the outer boundary is the 
same as two SLE$_{8/3}$ conditioned not to intersect. This shows that there might (maybe surprisingly) still exist a 
direct relation between $\mu$ and boundaries of some Brownian loops in this case. 

\medbreak

The measure $\mu_S$ defines (deterministic) global quantities related to the Riemann surface $S$. For instance, one can look at the 
total mass of the set of loops in a given homology class. 
These quantities therefore provide ways to describe the moduli space of Riemann surfaces.
Furthermore, the restriction property of $\mu$ makes it 
well-suited to the study of ``local deformations'' of the complex structure. 
In particular, suppose that two Riemann surfaces $S$ and $S'$ coincide almost everywhere, but are different ``inside 
a small circle $C$''. Then, the only difference between $\mu_S$ and $\mu_{S'}$ concerns just those loops that go through the domain  encircled by $C$. Hence, one can control the variation of the macroscopic observables defined by the measures $\mu$ with respect to ``infinitesimal variations'' of the Riemann structure.
We plan to investigate this further in upcoming work \cite {Wip}.

\medbreak

Recent attention (e.g. \cite {Zhan, BF}) has been devoted to the question of generalizing 
the definition of SLE to non-simply connected domains. More generally, one can try to define SLE 
on any Riemann surface. The present construction of $\mu$ on any Riemann surface together with the 
construction of conformal loop-ensembles \cite {Wls, ShW} (see also \cite {Wln3}) via Poissonian clouds 
of self-avoiding loops
gives a way to define families of SLE-loops on any Riemann surface (even if it does not have a boundary).
Note that in the case of surfaces with boundaries (and therefore finite Greens functions), a Poissonian cloud of 
Brownian loops leads to the same type of simple CLE, but for compact surfaces, the Brownian loop-clusters become dense and uninteresting, while the clusters of self-avoiding loops still make sense.
  
\section {Asymptotics for the function F}

\subsection {Motivation}
\label {section6.1}
A motivation for the present paper is the desire to understand better the 
behavior of long self-avoiding walks on planar lattices and the related conjectures. 
For some background and recent status on these conjectures, see \cite {LSWsaw}. 
Let us mention here the following conjecture from \cite {LSWsaw}: 
Consider a regular periodic planar lattice $L$ (to fix ideas, let us suppose that the lattice is
the square lattice, the triangular lattice or the hexagonal lattice). Recall that there exists a lattice-dependent
 constant 
$\lambda > 1$ 
such that the number of self-avoiding paths of length $n$ on $L$ that start at the origin behaves like
$\lambda^{n+o(1)}$ when $n \to \infty$.

Let us now define the discrete measure $\mu^{\delta}$ on the set of self-avoiding loops that one can draw on $\delta L$, that (for each $n$) puts a weight $n^{- \lambda}$ to each loop with $n$ steps.

\begin {conjecture}[\cite {LSWsaw}]
 $\mu^{\delta}$ converges when $\delta \to 0$ to a limiting measure on 
self-avoiding loops. 
\end {conjecture}

If this limiting measure is conformally invariant, then it is easy to see that it must be (a multiple of) 
our measure $\mu$.
So, if one believes in this conjecture and in conformal invariance in the scaling limit, 
then the fact that the definition 
of $\mu$ makes sense on any Riemann surface and satisfies restriction in a broad sense
should not be surprising.     

One may wonder if a further study of the (conjectural) limiting measure $\mu$ might be of some help to
understand the conjecture better. This is what motivates the next subsection.

\subsection {Asymptotics}

It would be nice to be able to have an explicit expression for the mass of the set of 
loops that go around an annulus in terms of its modulus (i.e. the explicit expression of $F$). This does not 
seem easy. 
We are now going to study the asymptotics of $F$ when $\rho \to 0$.

\begin {proposition}
For some constant $c_0$, $F(\rho)  \sim c_0 e^{-{5 \pi^2 }/ {4\rho } }$
when $\rho \to 0$.
\end {proposition}

In other words, if one considers a very thin tube of width $1$ around a very long smooth curve of length
$l$, the mass of the set of loops that go around the tube decays exponentially in $l$. This is of course reminiscent of the fact that for the discrete grid-model $\mu^\delta$, the mass decays exponentially in the number of steps.

As we shall see, the exponent $5 \pi^2 / 4$ is very closely related to the $5/8$ exponent for SLE$_{8/3}$ as derived in \cite {LSWr}. 

Before proving this fact, let us  first concentrate on 
the corresponding result for the bubble measure. For convenience, we 
will first use the unit disc rather than the upper half-plane.
Define $f(\rho)$ (up to a multiplicative constant) as the limit as $\eps \to 0$ of $\eps^{-2}$ 
times the probability that an 
SLE$_{8/3}$ from $1$ to $\exp ( i \eps)$ in $\U$ goes ``around'' (i.e. clockwise in this case) the 
disc $e^{-\rho} \U$ without hitting it. 
We want to control this probability uniformly with respect to $\rho$ when $\eps \to 0$.

\begin {lemma}
\label {l20}
For some positive constant $c$, $f(\rho) \sim 
c \rho^{-2} \exp ( - 5 \pi^2 / (4 \rho))$ when $\rho \to 0$.
\end {lemma}

\noindent
{\bf Proof of the lemma.}
It is convenient to use the Poissonian excursion representation of SLE$_{8/3}$ 
introduced in \cite {Wcrr}. Let us 
briefly recall how it goes in the present case of the unit disc. First, for each $\theta \in [0, 2 \pi)$, we
define the Brownian excursion measure ${\cal E}^\theta$ at $e^{i \theta}$ in $\U$ 
as the limit when $\delta \to 0$  
of $\delta^{-1} \times P^{BM}_{\exp (- \delta + i \theta)}$, where 
$P^{BM}_z$ denotes the law of a Brownian motion started from $z$ and stopped at its first exit of $\U$.
Then, 
define the excursion measure $ {\cal E} = \int_0^{2\pi} d \theta {\cal E}^\theta$. This is an infinite measure, supported on the set of Brownian paths in $\U$ that start and end on $\partial \U$ (each excursion has a starting point and an endpoint). It has been first defined in \cite {LW1} and for instance used in \cite {LW,Wcrr}. It possesses strong conformal invariance properties (see e.g. \cite {LW}). 

Then, take two points $u_1$ and $u_2$ on $\partial \U$, and let $\partial_1, \partial_2$ denote the two 
connected components of $\partial \U \setminus \{u_1, u_2\}$.
Consider a Poissonian realization of $c {\cal E}$ restricted to those excursions that have both their endpoints on 
$\partial_1$, where $c$ is a well-chosen constant. Consider the curve from $u_1$ to $u_2$ in $\U$ that delimits the boundary of the connected component (that contains $\partial_2$ as the other part of its boundary) of the complement in $\U$ of the union of all these excursions in $\U$. 
Then (see \cite {Wcrr}), for this well-chosen definition of $c$, the curve $\eta$ is exactly an SLE$_{8/3}$ from $u_1$ to $u_2$. 

In the case where $u_1 =1$ and $u_2 = e^{i \eps}$ and $\partial_1$ is the ``long arc'' 
$\{ e^{i \theta}, \theta \in [\eps, 2 \pi] \}$, then clearly, if we want $\eta$ to go around the disc 
$e^{- \rho} \U$, necessarily, none of the excursions is allowed to hit $e^{- \rho} \U$.   

The ${\cal E}^\theta$ mass of the set of excursions that hit $e^{-\rho} \U$ is equal to 
$1/\rho$ (just because $\log |Z|$ is harmonic for a planar Brownian motion $Z$). Hence, the 
${\cal E}$-mass of the set of excursions that hit $e^{-\rho } \U$  is equal to 
$2 \pi / \rho$. Also, the ${\cal E}$-mass of the set of excursions that hit $e^{- \rho} \U$  and 
has their starting point (resp. end point) in $\partial_1$ is $(2 \pi - \eps)/ \rho$. 
Hence, the ${\cal E}$-mass of the 
set of excursions that hit $e^{-\rho} \U$ and have both their endpoints in $\partial_1$ is greater than 
$(2 \pi - 2 \eps) / \rho$.
It follows that the probability $p_{\eps, \rho}$ that none of the excursions used to 
construct the SLE from $u_1$ to $u_2$ in $\U$
does hit $e^{- \rho} \U$ satisfies
$$ 
\exp \{ - 2 \pi c / \rho \} \le p_{\eps, \rho} \le \exp \{-2 (\pi - \eps) c / \rho \}
$$     
(recall that for a Poissonian sample of intensity measure $m$, the probability that in the Poissonian  
sample, no item belongs to the set $A$ is $\exp \{- m (A)\}$).

We say that an excursion in the annulus $\U \setminus e^{-\rho} \U$ goes the easy way 
if it disconnects $1$ from $e^{-\rho} \U$ in $\U$. Otherwise, we say that it 
goes the ``hard way''.
When no excursions intersect the disc $e^{- \rho} \U$, then the SLE itself does not intersect the smaller disc either.
It then remains to estimate the (conditional) 
probability that the SLE goes ``the hard way'' (i.e. ``around $e^{-\rho} \U$'').
This means exactly that none of the two following events happens:
\begin {itemize}
\item
There exists a Brownian excursion that goes the easy way.
\item
There exists two Brownian excursions that do not go the easy way, such that their union 
 disconnects $1$ from $e^{-\delta} \U$
in $\U$.
\end {itemize}
We know already that this conditional probability decays like $g(\rho) \eps^2$ when $\rho$ is 
 fixed and $\eps \to 0$, for some positive $g(\rho)$ because of the convergence of the 
 renormalized SLE measure to the bubble measure (i.e. $g(\rho) \exp (- 2  \pi c / \rho)$ is the mass for the 
 bubble measure of the set of bubbles that go around $e^{-\rho} \U$).
 
 On the other hand, it is easy to note that because of scaling, 
 when $\eps \to 0$, the probabilities of this event for $\eps, \delta$ and the event for $2\eps, 2 \delta$ become
 very close (for instance by comparing this event with the event 
 in the strip $\R \times [0, \delta])$). Hence, we get that $g(\rho)= cst / \rho^2$.
 
 We have now proved that for some constants $c_1, c_2$,
 $$f(\rho) \sim c_1 \rho^{-2} \exp \{ - c_2 / \rho \} $$
 when $\rho \to 0$.
 In fact, it is possible to identify what the constant $c_2$ is because we know that the 
chordal restriction exponent of SLE$_{8/3}$ is $5/8$ (see \cite {LSWr}). In particular, this enables for instance to 
estimate the probability that an SLE from $-1$ to $1$ in the 
unit disc stays in $\U \setminus (e^{-\rho}\U \cup [-i,0])$ and the leading order term when $\rho \to 0$ 
can be interpreted in terms of Brownian excursions in $\U$. 
The upshot is that $c_2 = 5 \pi^2 / 4$.
Intuitively, going around the annulus is at first order like crossing the 
rectangle $[0, 2 \pi] \times [0, \rho]$. By scaling, this is like 
crossing the rectangle $[0,  L=2 \pi^2 / \rho ] \times [0, \pi]$, and for an SLE$_{8/3}$, this has probability of
order $\exp ( - 5L/8)$.
\qed

\medbreak
\noindent
{\bf Proof of the proposition.}
It now remains to deduce the result for the loop measure $\mu$ from the 
result on bubbles.
For this, we are going to use the decomposition of $\mu$ along a Loewner chain.
Suppose that $\eta : [0, T] \to \overline \U$ now denotes a continuous simple path in the unit disc
$\U$, started from $\eta_0=1$. We asume that it is parametrized from ``capacity'' at the origin as radial
Loewner chains usually are i.e. that if we define the conformal maps $\varphi_t$ from 
$U_t=\U \setminus \eta [0,t]$ onto $\U$ such that $\varphi_t (0) = 0$ and
$\varphi (\eta_t) = 1$, then $| \varphi_t'(0)| = \exp t$.
In fact, we may want to think of $\eta= \eta^\eps$ as a spiral that winds very fast around the origin 
(say a time-reparametrization of $\exp (- t + i t / \eps )$ for a very small $\eps$).
Then, for small $\eps$, the domain $U_t$ is (conformally speaking) very close to the disc $e^{-t} \U$.

The function $f(r)$ can be viewed as the $bubb$-measure of the set of bubbles that go around the 
disc $e^{-\rho}$, where $bubb$ is ``the bubble measure in the unit disc'', started at the boundary point $1$.
Exactly as in the case of the upper-half plane, one can prove that 
$$
\mu_\U 1_{\gamma \cap \eta [0,T] }  = \int_0^T dt \varphi_t^{-1} \circ bubb .
$$
Hence, it follows that 
$$
\mu_\U ( \gamma \cap \eta [0,T] \not= \emptyset \hbox { and } \gamma
\hbox { surrounds } e^{-\rho} \U ) 
 = \int_0^T f( r(t)) dt 
$$
where $r(t) = \rho ( \U \setminus (\eta [0,t] \cup e^{-\rho} \U)) $.
If we take $T= \rho$, and $\eta=\eta^\eps$ as before, and then let 
$\eps \to 0$, we get immediately that 
$$
F( \rho) = \mu_\U  (\gamma
\hbox { surrounds } e^{-\rho} \U ) 
=
\int_0^\rho f(r) dr.
$$
In particular, when $\rho \to 0$,
$$ 
F(\rho) \sim c_1 \exp ( - 5 \pi^2 / 4 \rho ).
$$
 \qed

\section {Relation to critical percolation}

We are going to conclude this paper by briefly mentioning in a more informal style
the direct relation between $\mu$ and the measures 
on outer perimeters of critical percolation clusters in their scaling limit.
This will show that, as in the case of Brownian loops, the inside perimeters and the 
outside perimeters of critical percolation clusters have the same ``law'', and also that 
outside perimeters are invariant (in law) with respect to perturbations of the complex 
structure ``inside''.

Let us consider critical percolation on the triangular lattice in the entire plane (we restrict ourselves to this 
lattice since it is at present the only one for which conformal invariance is proved).
We can view it as a random coloring of the faces of a hexagonal grid, where each face is colored 
independently in black or white with probability $1/2$. 
One way to encode the obtained configuration is to define the family of  loops on the 
hexagonal grid that separate white connected components and black connected components. 
These are the interfaces for this model. Each of the loops is self-avoiding on the scale of the lattice and 
the family of all loops describes fully the percolation configuration.

\begin{figure}[t]
\begin{center}
\includegraphics [width=3.5in]{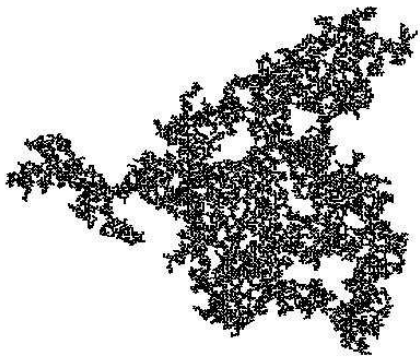}
\caption{A percolation cluster (low resolution)}
\end{center}
\end{figure}

For convenience, we are going to introduce a small $\delta$, and look at the percolation configuration 
on a grid of mesh-size $\delta$ (we shall then let $\delta \to 0$).
It is easy to see using the Russo-Seymour-Welsh theory that for any bounded domain $D$ and any 
$\eps >0$, the mean number of loops included in $D$ and of diameter greater than $\eps$ 
remains bounded and bounded away from $0$ (except of course if $\eps$ is greater than the 
diameter of $D$...) as $\delta \to 0$.

Furthermore, it is also easy to see that the loops will ``not remain self-avoiding'' on the large scale, i.e., 
that on a large loop, one can find neighboring sites $x$ and $x'$ on the hexagonal lattice of mesh-size $\delta$
 and on the same 
loop such that the loop goes (twice) to a macroscopic distance between its visits of $x$ and $x'$.

Smirnov \cite {Sm} proved  that in the scaling limit crossing probabilities 
of conformal rectangles for this percolation model converge to a (conformally invariant) limit
predicted by Cardy \cite {Ca} when $\delta$ goes to zero.
Camia and Newman \cite {CN} showed that, using Smirnov's result, how one can indeed deduce 
convergence when $\delta \to 0$ of the joint law of all the discrete loops, towards a continuous 
family of loops, that are also  described in terms of the
SLE$_6$ processes that had been studied in \cite {LSW1,LSW2}.
This continuous family of loops in some sense describe the joint law of the clusters in the scaling limit.
In particular, one loop out of two is the outer boundary of a white cluster, the other ones being ``inside boundaries''
of white clusters.
Let us call $(o_j, j \in J)$ these outer boundaries.

\begin{figure}
\begin{center}
\includegraphics [width=3.5in]{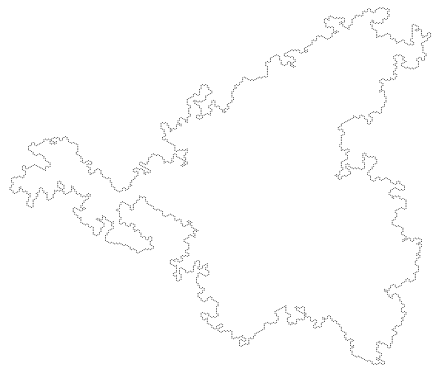}
\caption{Its outer loop}
\label{clusterloop}
\end{center}
\end{figure}

These continuous loops are not self-avoiding
 (in fact, they have many double points: the dimension of each loop is a.s. $7/4$ 
while the dimension of their ``boundaries'' is $4/3$, see e.g. \cite {Be1}). Let us call $O_j$ the self-avoiding loop obtained by taking the ``outer boundary'' (i.e. the boundary of the unbounded connected component of the 
complement) of $o_j$. The conformal invariance statement can then be phrased as follows:

For each simply connected domain $D$, define the law joint law $P_D$ of the set of loops 
$(o_j, j \in J_D)$ that do stay in $D$. Then, if $\Phi$ is a conformal map from $D$ onto another simply connected domain $\tilde D$, then 
the law of $(\Phi (o_j), j \in J_D)$ is identical to $P_{\tilde D}$.

Clearly, this implies also the same statement for the outer boundaries $O_j$ instead of the boundaries $o_j$.
Hence, the probability measures $P_D$ satisfy a type of conformal restriction property. In fact, if we now define 
for each measurable set $A$ of self-avoiding loops
$$
\pi ( A ) =  E ( \# \{ j \in J \ : \ O_j \in A \} )
,$$
then $\pi$ is a measure on the set of self-avoiding loops in the plane, that does satisfy conformal restriction.
Hence, it is equal to a constant times the measure $\mu$ that we have defined in this paper.
In other words, the
 measure $\mu$ can also be constructed as the outer boundary of the scaling limits of 
percolation cluster boundaries. 
The multiplicative constant $c$ does here not seem to be easy to determine.

Note that in fact, as we are focusing only on these outer boundaries $O_j$, we do not really need the convergence of 
the full discrete loops to the full loops $o_j$. Indeed, it is not very hard using some a priori estimates for
percolation (i.e. the 5-arm estimate) to see that the $O_j$'s are the limits when $\delta\to 0$ of the 
discrete outer perimeters of percolation clusters, where the outer perimeter of a cluster is the outer boundary of 
the domain obtained by ``filling in all fjords of width one''. More precisely, consider the unbounded connected component of the complement of a cluster. 
We then just keep the subdomain of this component consisting of points $x$ such that there exist two disjoint paths 
from $x$ to infinity that avoid the cluster. The boundary of this subdomain is the discrete outer perimeter (these discrete curves had been studied and simulated by physicists before, see e.g. \cite {ADA} and the references therein).
Then, the law of the joint family of discrete outer perimeters
$(O_j^\delta, j \in J^\delta)$ converges when $\delta \to 0$ to that of $(O_j, j \in J)$, so  that
for some constant $c$ and any (reasonable) set of (macroscopic but bounded) loops (for instance the set of loops that go around the hole in a given conformal annulus), 
$$ 
\mu (A) = c \lim_{\delta\to 0} E ( \# \{ j \in J^\delta \ : \ O_j^\delta \in A \} ) 
.$$     
The shape of a percolation cluster outer perimeter has therefore exactly the same law as that of the 
outer boundary of a Brownian loop. Note again that the proof of this result relies solely on the 
conformal invariance of both models and on their restriction property.

\begin{figure}[t]
\begin{center}
\includegraphics [width=3.5in]{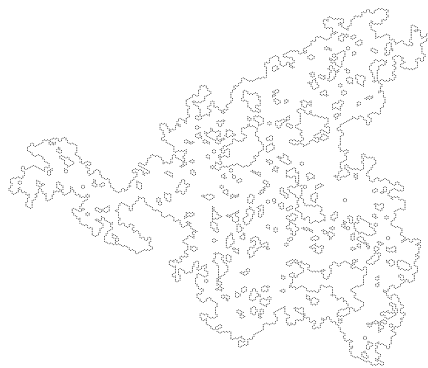}
\caption{Its outer and inner loops}
\label{clusterloops}
\end{center}
\end{figure}

The definition of the discrete outer perimeter is very non-symmetric, and
 it is non-local in the discrete case: Because of the ``closure of fjords procedure'', 
 it does not suffice to look at the status of the hexagons  
that touch a given self-avoiding loop to decide if it is an outer perimeter. One needs to make sure that the ``fjords do not meet in the inside''.
Hence, the inside/outside symmetry statement for $\mu$ is not obvious at all from this perspective either
(and this statement is based on SLE considerations).
In fact, it has rather surprising consequences, just as for the Brownian loop approach.
It says that outer perimeters of percolation are (in law)
not affected by change of the conformal structure in their inside.

As for the Brownian measure $M$ in Corollary \ref {c15},
 one can also characterize  $\mu$ in terms of all perimeters (not only the outer ones)
and this changes the constant $c$ by a factor of $2$. We safely leave this to the reader. 
It is not difficult to check that the joint law (i.e. measure) of the boundaries for one percolation cluster
is not the same as the joint law (i.e. measure) of the boundaries for one Brownian loop. For instance, the existence 
of cut-points for Brownian motions shows that in a Brownian loop, there exist inner boundaries that intersect 
the outer boundary, and that the dimension of these intersection points is a.s. 3/4 (see \cite{LSW2}).
 For the scaling limit of percolation clusters, this is not the case (it would correspond to a 
``6 arms'' exponent -- recall that we are looking at outer boundaries). It is just the shape of ``one prescribed''
boundary that coincide for percolation clusters and for Brownian loops.

Finally, note that in the case of Riemann surfaces of higher genus, the cut-point issue raised in the Brownian loop case (i.e. the fact that the outer boundary can have double points) does not occur
in the case of outer boundaries of percolation clusters, because the perimeters' sides have a definite ``color''.

\bigbreak

{\bf Acknowledgements.} I would like to thank Pierre Nolin for the pictures of the Brownian loop, the percolation cluster and their boundaries.

\medbreak

{{Universit\'e Paris-Sud\\
Laboratoire de Math\'ematiques, Universit\'e Paris-Sud
\\
B\^at. 425, 91405 Orsay cedex, France \\}}

\medbreak
 
{DMA,
 Ecole Normale Sup\'erieure, \\
45 rue d'Ulm,\\
75230 Paris cedex, France}\\

\medbreak
email: wendelin.werner@math.u-psud.fr

\medbreak

{Work supported by the Institut Universitaire de France}

\end {document}